\pdfoutput=1
\documentclass{amsart}
\usepackage{colortbl}
\usepackage{amsmath,amsthm,amssymb}
\usepackage{graphicx}
\usepackage{harvard}

\renewcommand{\L}{\ensuremath{\mathcal{L}}}

\newtheorem{definition}{Definition}
\newtheorem{lemma}[definition]{Lemma}
\newtheorem{theorem}[definition]{Theorem}
\newtheorem{corollary}[definition]{Corollary}

\newtheorem{remark}[definition]{Remark}
\newtheorem{proposition}[definition]{Remark}

\numberwithin{equation}{section}

\renewcommand{\(}{\left(}
\renewcommand{\)}{\right)}
\renewcommand{\[}{\left[}
\renewcommand{\]}{\right]}

\def\F{{\mathcal F}}
\def\L{{\mathcal L}}

\def\G{\tilde{\mathcal L}}

\def\R{\mathbb R}
\def\C{\mathbb C}
\def\N{\mathbb N}

\def\Z{\mathbb Z}

\def\P{\mathrm P}

\def\a{\bar a}
\def\b{\bar b}
\def\c{\bar c}

\def\partI#1#2{\frac{\partial#1}{\partial#2}}
\def\partII#1#2{\frac{\partial^2#1}{\partial#2^2}}

\def\espcond#1#2#3{E^{#1}\left[\begin{array}{l|r}#2&#3\end{array}\right]}

\begin{document}

\title{Laplace Transforms for Integrals of Markov Processes}

\author{Claudio Albanese}
\address{Department of Mathematics, Imperial
College, London, U.K.}

\author{Stephan Lawi}
\address{Laboratoire de
Probabilit\'es et Mod\`eles Al\'eatoires, CNRS (UMR 7599), Paris,
France}

\thanks{Supported in part by the Natural Science
and Engineering Council of Canada under grants RGPIN-171149.}

\date{\today}
\maketitle

\begin{abstract}

Laplace transforms for integrals of stochastic processes have been
known in analytically closed form for just a handful of Markov
processes: namely, the Ornstein-Uhlenbeck, the Cox-Ingerssol-Ross
(CIR) process and the exponential of Brownian motion. In virtue of
their analytical tractability, these processes are extensively used
in modelling applications. In this paper, we construct broad
extensions of these process classes. We show how the known models
fit into a classification scheme for diffusion processes for which
Laplace transforms for integrals of the diffusion processes and
transitional probability densities can be evaluated as integrals of
hypergeometric functions against the spectral measure for certain
self-adjoint operators. We also extend this scheme to a class of
finite-state Markov processes related to hypergeometric polynomials
in the discrete series of the Askey classification tree.

\end{abstract}

\section{Introduction}
Let $(X_t)_{t\ge0}$ be a time-homogenous, real-valued Markov process
on the filtered probability space $(\Omega, \{\F_t\}_{t\ge0}, {\rm
P})$ and consider the Laplace transform $L_{T-t}(X_t, \vartheta)$
defined as follows:
\begin{equation}
L_{T-t}(X_t,\vartheta) = \espcond{\P}{e^{- \vartheta \int_t^T
\phi(X_s) ds}\ q(X_T)}{\F_t}
\end{equation}
where  $t\le T$, $\vartheta\in\C$ and $\phi, q:\R\to\R$ two Borel
functions. In this paper, we address the question of whether it is
possible to compute the Laplace transform $L_{T-t}(X_t, \vartheta)$
in analytically closed form. Our work builds upon several streams of
research often motivated by applications to various fields of
Physics and Finance, and unifies them to obtain a broad
classification scheme for Laplace transforms expressible in analytic
closed form.\\

For $q\equiv1$, a class of examples for which analytic closed form
solutions are available is represented by the so called {\it affine
models} which are characterized by a representation of the form
\begin{equation}
L_{T-t}(X_t, \vartheta) = e^{m( T-t, \vartheta) X_t  + n(T-t,
\vartheta)}.
\end{equation}
The archetypical affine models are based on diffusion processes and
are described by stochastic differential equations of the form
\begin{equation}\label{eq:sde_affine}
d X_t = (a - b X_t) dt + \sigma X_t^\beta dW_t
\end{equation}
where $a, b, \sigma$ are constants and $\beta = 0$ or $\frac{1}{2}$.
The case $\beta = 0$ corresponds to the Gaussian Ornstein-Uhlenbeck
process \cite{Va1977} and the case $\beta = \frac{1}{2}$ corresponds
to the Cox-Ingersoll-Ross (CIR) process \cite{CIR1985}. The case of
the CIR process was generalized to bridges by Pitman and Yor in
\cite{PY1982}. It has been shown in \cite{BS1994a}, \cite{BS1994b}
and in \cite{DG2003} that any affine process which is a
time-homogenous, nonnegative diffusion is necessarily of the CIR
type. However, there are also affine processes with jumps. General
non-negative affine processes correspond to the so called
conservative CBI-processes (continuous state branching processes
with immigration) and have been well studied, among others, by
Kawazu and Watanabe in \cite{KW1971} and Filipovi\'c in
\cite{Fi2001a}.

An extension of the affine class, known as the quadratic class,
postulates the Laplace transform being of the form
\begin{equation}
L_{T-t}(X_t, \vartheta) = e^{l(T-t, \vartheta) X_t^2 + m( T-t,
\vartheta) X_t  + n(T-t, \vartheta)}.
\end{equation}
The first examples of quadratic models appeared in the double square
root model of Longstaff in \cite{Lo1989} and in the nonlinear
equilibrium model by Beaglehole and Tenney in \cite{BT1992}. Rogers
\cite{Ro1997} also uses examples where the pricing kernel is a
quadratic function of the Markov process. Most recently, Filipovi\'c
\cite{Fi2001b} proved that if one represents the forward rate as a
polynomial function of the diffusion process, the maximal consistent
order of the polynomial is two. Consistency in this context means
that the interest rate model will produce forward rate curves
belonging to the parameterized family. Finally, Leippold and Wu
\cite{LW2002} formulated a general asset and derivative pricing
framework for the quadratic class.

A separate class of models for which the Laplace transform can be
expressed in analytically closed form is represented by the
exponential Brownian motion of equation
\begin{equation}\label{eq:sde_geob}
d X_t = \mu  X_t dt + \sigma  X_t dW_t
\end{equation}
where $\mu, \sigma >0$ are positive constants. This case was first
considered by Yor in \cite{Yo1992} who arrived to an expression
involving a triple integral. An earlier related result for bond
prices given in terms of an integral over modified Bessel functions
was formulated by Dothan in \cite{Do1978}. As an alternative, Geman
and Yor \cite{GY1993} derive a closed-form expression for the
Laplace transform in terms of confluent hypergeometric functions
(see Donati-Martin et al. \cite{DMGY2001} and Yor \cite{Yo2001} for
further references). For applications to finance, one needs to
compute the inverse Laplace transform, for which numerical methods
have been developed by Geman and Eydeland \cite{GE1995}, Fu et al.
\cite{FMW1998}, Craddock et al. \cite{CHP2000}, Shaw \cite{Sh2002}.
Dufresne in \cite{Du2000} and Linetsky in \cite{Li2003} develop
analytical methods and alternative expansions.\\

In this article, we obtain far reaching extensions of the
representation formula for the Laplace transform for the integrals
of stochastic processes over geometric Brownian motions. The key
idea is to seek expansions of similar form as those in \cite{Du2000}
and \cite{Li2003} but expressed in terms of more general
hypergeometric functions, see \cite{AS1972}, as pioneered by Wong in
\cite{Wo1964}. Let us recall that general hypergeometric functions
are denoted as follows:
\begin{equation}
_pF_q(\alpha_1,\ldots,\alpha_p;\gamma_1,\ldots,\gamma_q;z)
\end{equation}
for $p\le q+1, \gamma_j \in \C\setminus -\Z_+$, and are represented
by the following Taylor expansion around $z=0$:
\begin{equation}
_pF_q(\alpha_1,\ldots,\alpha_p;\gamma_1,\ldots,\gamma_q;z) =
\sum_{n=0}^\infty\ \frac{(\alpha_1)_n \ldots (\alpha_p)_n}{
(\gamma_1)_n \ldots (\gamma_q)_n}\ \frac{z^n}{n!}.
\end{equation}
The Kummer functions in the work by Geman and Yor \cite{GY1993} are
in the family $_1F_1$ of the so called confluent hypergeometric
functions. Gaussian hypergeometric functions are in the family
$_2F_1$ and admit the functions of type $_1F_1$ as limits. Both
Gaussian and confluent hypergeometric functions solve differential
equations of the Fuchsian class, see \cite{Gr1984}. Specifically
\begin{equation}\label{eq:hypg}
z(1-z) \ _2F_1''(\alpha, \beta; \gamma; z) + (\gamma - (1+\alpha +
\beta) z) \ _2F_1'(\alpha, \beta; \gamma; z) - \alpha \beta \
_2F_1(\alpha, \beta; \gamma; z) = 0
\end{equation}
and
\begin{equation}\label{eq:hypc}
z  \ _1F_1''(\alpha; \gamma; z) + (\gamma - z) \ { }_1F_1'(\alpha;
\gamma; z) - \alpha \ _1F_1(\alpha; \gamma; z) = 0.
\end{equation}
In general, higher order hypergeometric functions are not associated
to a differential equation. However, in some particularly important
cases, they provide solutions of finite difference equations. The
Askey classification scheme, see \cite{KS1998} and \cite{Sz1959},
gives a complete list of all orthogonal polynomials solving either a
differential or a finite difference equation and in addition satisfy
a recurrence relation. All of these polynomials descend as
particular or limiting cases from the so-called Racah polynomials,
which are particular cases of the hypergeometric functions $_4F_3$.

We first consider the case of diffusion processes and next the case
of finite state Markov processes. In the diffusion case, we
construct a classification scheme based on reduction to eigenvalue
problems admitting solutions within the class of Gaussian and
confluent hypergeometric functions $_2F_1$ and $_1F_1$. In the
second case, the problem is more difficult for several reasons, as
there is no discrete equivalent of a theory of Fuchs type equations
and, in addition, the groups of conformal transformations and
diffeomorphisms do not extend to lattices. What we do in the
discrete case is to take the moves from the Askey classification
scheme for orthogonal polynomials and show how to extend the
previous spectral decomposition to the case of Meixner, dual Hahn
and Racah polynomials, which are special cases of $_2F_1$, $_3F_2$
and $_4F_3$ hypergeometric functions.\\

For a diffusion process, on a domain $D_x\subset\R$, of the form
\begin{equation}\label{eq:sde}
d X_t = \mu(X_t) dt + \sigma(X_t) dW_t,
\end{equation}
we define the transitional probability density $p_{T-t}(x, y)$ as
the density of the Markov semigroup of the process $X_t$:
\begin{equation}
\espcond{\P}{f(X_T)}{\F_t} = \int_D f(y)\ p_{T-t}(x, y) dy.
\end{equation}
We are interested in building a classification scheme for the drift
and volatility functions $\mu(x), \sigma(x)$ such that the
calculation of both functions $p$ and $L$ can be reduced to
computing an integral over hypergeometric functions. The
transitional probability density and Laplace transform can be
computed in terms of the spectral resolution for the infinitesimal
generator of the process $X_t$
\begin{equation}\label{eq:L}
\L = \frac{\sigma(x)^2}{2} \frac{\partial^2}{\partial x^2} + \mu(x)
\frac{\partial}{\partial x}
\end{equation}
and the Feynman-Kac operator
\begin{equation}\label{eq:G}
\G = \L - \vartheta \phi(x) = \frac{\sigma(x)^2}{2}
\frac{\partial^2}{
\partial x^2} + \mu(x) \frac{\partial}{\partial x} - \vartheta
\phi(x).
\end{equation}
In fact, we have that
\begin{equation}\label{eq:probk}
p_{T-t}(x, y) = e^{(T-t) \L}(x, y)
\end{equation}
and
\begin{equation}\label{eq:Genf}
L_{T-t}(x, \vartheta) = \int_D q(y)\ e^{(T-t)\G}(x, y) dy.
\end{equation}
As we show in detail in Section \ref{sec:Diffusion}, these operators
are conjugated by a non-singular transformation to self-adjoint
operators which admit a spectral resolution. The calculation of the
transitional probability density in (\ref{eq:probk}) and the Laplace
transform in (\ref{eq:Genf}) is thus reduced to the resolution of
the differential eigenvalue problems
\begin{equation}\label{eq:eigp}
\L\ f(x) = \lambda\ f(x)\quad, \quad\G\ \bar f(x) = \bar \lambda\
\bar f(x).
\end{equation}
To properly define the classification scheme, we specify by what
means the reduction can be accomplished.

\begin{definition}\label{def:3op}
The problem of finding the transitional probability density and the
Laplace transform for the process in (\ref{eq:sde}) is said to be
{\rm reducible to a spectral integral over hypergeometric functions}
if the two eigenvalue problems in (\ref{eq:eigp}) can be recast in
the form of a differential equation for hypergeometric functions
such as either (\ref{eq:hypc}) or (\ref{eq:hypg}), by means of a
combination of the following three operations $T_i$:
\begin{enumerate}
\item $\ T_Z$: change of variable $x \mapsto z = Z(x)$
where $Z(x)$ is a diffeomorphism $Z:D_x\to D_z$ such that
\begin{equation}
\L_x\mapsto \L_z\quad,\quad \G_x\mapsto \G_z,
\end{equation}
\item $\ T_h$: gauge transformation associated to a strictly positive
function $h$ such that
\begin{equation}
\L \mapsto h^{-1} \L\ h\quad , \quad\G\mapsto h^{-1}\G\ h,
\end{equation}
\item $\ T_{\gamma^2}$: left-multiplication by a strictly
positive function $\gamma^2$ such that
\begin{equation}
\L \mapsto \gamma^2 \L \quad , \quad\G\mapsto \gamma^2\G.
\end{equation}
\end{enumerate}
\end{definition}
The third kind of transformations was first recognized in its
generality by Natanzon in the article \cite{Na1971} on integrable
Schr\"odinger equations, see also Milson's paper \cite{Mi1998}. The
following theorem gives a concise statement of our main
classification result:

\begin{theorem} [First Classification Theorem] \label{thm:Diffusions}
The most general reducible diffusion process (up to diffeomorphism)
according to Definition \ref{def:3op} can be constructed as follows:
\begin{enumerate}
\item four second order polynomials in $x$: $A(x),\ Q(x,\vartheta),\
R(x),\ S(x)$, such that $A(x)$ belongs to the set
$\{1,x,x(1-x),x^2+1\}$ and $R(x)\ge0$;

\item conditions for the stochastic process $X_t$
on the boundary of the domain $D_x$, specifying the relative
probability of reflection versus absorption upon hitting the
boundary;

\item a solution of the following equation in $D_x$ for some
$\xi\in\R$:
\begin{equation}\label{eq:T1v}
\frac{A(x)^2}{R(x)} h''(x) + \frac{S(x)}{R(x)} h(x) = \xi h(x)
\end{equation}
The function $h(x)$ is a linear combination of hypergeometric
functions of the confluent type $_1F_1$ if $A(x)\in \{1,x\}$ and of
the Gaussian type $_2F_1$ if $A(x)\in \{x(1-x), x^2+1\}$.
\end{enumerate}
The process associated to this choice is given by a generic solution
to the following stochastic differential equation on the domain
$D_x$:
\begin{equation}\label{eq:dx}
dX_t = 2\frac{h'(X_t)}{h(X_t)}\ \frac{A(X_t)^2}{R(X_t)}\ dt +
\frac{\sqrt{2}A(X_t)}{\sqrt{R(X_t)}}\ dW_t,
\end{equation}
with the boundary conditions above. The Laplace transform is
specified by
\begin{equation}\label{eq:phi}
\phi(x) = \frac{Q(x,\vartheta) }{\vartheta R(x)}.
\end{equation}
\end{theorem}
\vspace{0.5cm}

The proof of this theorem is in Section \ref{sec:Diffusion}. These
constructs are based on spectral analysis techniques for which we
refer to the book by Reed and Simon \cite{RS1975b}. We also make the
spectral analysis more explicit and list the expressions for the
transitional probability density and Laplace transforms in terms of
the kernel of semigroups generated by integrable quantum
Schr\"odinger operators. In Section \ref{sec:Examples}, we
specialize further and re-discover the known cases of processes
built upon the geometric Brownian motion and on the
Ornstein-Uhlenbeck and CIR processes, along with some interesting
extensions.

In Sections \ref{sec:ContPoly} and \ref{sec:DiscPoly}, we restrict
the framework to the special case of hypergeometric polynomials. In
the discrete case, studied in Section \ref{sec:DiscPoly}, the
process $X_t$ takes on only a discrete set of values, as opposed to
following a diffusion process. In this class of models, we base our
analysis on the Askey-Wilson theory of orthogonal polynomials, see
\cite{AW1979} and \cite{Sz1959}. We briefly review the basic
notions, following \cite{KS1998}.

\begin{definition}
An {\rm orthogonal system of polynomials} is given by a sequence of
polynomials $Q_n(x)$ of order $n$ for $n \in\N$ on the interval
$D\subseteq\R$ which satisfies an orthogonality condition of the
form
\begin{equation}\label{eq:Ortho}
\int_D Q_n(x) Q_m(x) \rho(dx) = d_n^2 \delta_{nm}\quad, \quad
n,m\in\N,
\end{equation}
where the $d_n$ are constants and $\rho(dx)$ is a given measure. One
distinguishes between {\rm continuous polynomials} whereby
$\rho(dx)$ is absolutely continuous with respect to the Lebesgue
measure, i.e.
\begin{equation}
\rho(dx) = w(x) dx
\end{equation}
for some weight function $w(x)$, and {\rm discrete polynomials} for
which
\begin{equation}
\rho(dx) = \sum_{i=0}^N w(x_i) \delta(x - i)\quad, \quad N\in\N.
\end{equation}
\end{definition}
All orthogonal polynomials satisfy a {\it three-term recurrence
relation} of the form
\begin{equation}
x Q_n(x) = A_n Q_{n+1}(x) - B_n Q_n(x) + C_n Q_{n-1}(x)
\end{equation}
where $n\ge1$, $A_n>0$, $C_n\ge0$ and $B_n\in\R$. Together with the
conditions $Q_{-1}(x) = 0$ and $Q_0(x)=1$, all the $Q_n(x)$ can be
determined based on this recurrence relation. The converse is also
true and is known as the Favard theorem, see \cite{Ch1976}.
Moreover, they satisfy the following eigenvalue equation,
\begin{equation}\label{eq:EV}
\L Q_n(x) = \lambda_n Q_n(x),
\end{equation}
for $\L$ a second-order differential operator in the continuous case
or a finite difference operator in the discrete case.

The reducibility condition in Definition \ref{def:3op} is mirrored
by the following (inequivalent) one which refers to orthogonal
polynomials in the continuous series as opposed to Gaussian
hypergeometric functions:
\begin{definition}\label{def:OP}
The problem of finding the transitional probability density in
(\ref{eq:probk}) and the Laplace transform in (\ref{eq:Genf}) is
said to be {\rm reducible to a spectral integral over orthogonal
polynomials} if the two eigenvalue problems in (\ref{eq:eigp}) have
the same orthogonal polynomials as eigenfunctions.
\end{definition}

We first restrict the framework to continuous orthogonal polynomials
that have as generator a second-order differential operator,
\begin{equation}\label{eq:ContL}
\L = \frac{\sigma^2}{2} A(x)\partII{}{x} +(a-bx)\partI{}{x},
\end{equation}
acting on the Hilbert space $L^2(D,\rho)$, where
$A(x)\in\{1,x,x(1-x)\}$, $a\in\R$ and $b,\sigma>0$. This class
consists of the Hermite, Laguerre and Jacobi polynomials up to
diffeomorphism. Our main result concerning the latter continuous
orthogonal polynomials can be stated as follows:

\begin{theorem}[Second Classification Theorem]\label{thm:ContPoly}
The most general reducible diffusion process (up to diffeomorphism)
in the sense of Definition \ref{def:OP}, has infinitesimal generator
$\L$ given by (\ref{eq:ContL}). Its transitional probability density
can be expressed as
\begin{equation}\label{eq:ContPolyU}
p_{T-t}(x, y) = \sum_{n=0}^\infty \frac{e^{\lambda_n(T-t)}}{d_n^2}
Q_n(x;a,b) Q_n(y;a,b) w(y).
\end{equation}
Furthermore, for some parameters $\a,\b,C\in\R$,
\begin{equation}\label{eq:ContPolyphi}
\phi(x) = C + \big((a-bx)-(\a-\b x)\big)\frac{ A'(x)}{2\vartheta
A(x)}+ \frac{(\a-\b x)^2-(a-bx)^2}{2\vartheta\sigma^2 A(x)}
\end{equation}
and the Laplace transform is given by the following convergent
series:
\begin{equation}\label{eq:ContPolyG}
L_{T-t}(x, \vartheta) =\exp\(\int^{x} \frac{(\a-\b y)-(a-b
y)}{\sigma^2  A(y)}dy\)\ \sum_{n=0}^\infty e^{\bar\lambda_n(T-t)}
z_n Q_n(x;\a,\b).
\end{equation}
The coefficients $z_n$ are given by:
\begin{equation}\label{eq:ContPolyzn}
z_n = \frac{1}{\bar d_n^2}\int_{D_x} q(x)\exp\(-\int^x \frac{(\a-\b
y)-(a-b y)}{\sigma^2 A(y)}dy\) Q_n(x;\a,\b)\bar\rho(dx).
\end{equation}
\end{theorem}
\vspace{0.5cm}

In Section \ref{sec:ContPoly} we present the proof of this
alternative classification scheme based on orthogonal polynomials of
the continuous series. This discussion sets the premise for the
extension of the result to orthogonal polynomials in the discrete
series. Discrete orthogonal polynomials are characterized by a
finite difference generator on the Hilbert space $l^2(\Lambda_N,w)$
with $\Lambda_N$ the set $\{0,\ldots,N\}$.
\begin{definition}
Let $\Delta^h$ and $\nabla_+^h$ denote the difference operators
defined as follows:
\begin{equation}
\Delta^h y(x) = y(x+h)-2y(x)+y(x-h)\quad,\quad\nabla_+^h y(x) =
y(x+h)-y(x).
\end{equation}
\end{definition}
Using these operators, the finite difference generator $\L$ takes
the form
\begin{equation}\label{eq:DiscL}
\L = -D(x)\Delta^1 + \big(D(x)-B(x)\big)\nabla_+^1,
\end{equation}
where $B(x)$ and $D(x)$ are rational functions of at most fourth
order in the numerator and at most second order in the denominator.
Our main result in the discrete case can be stated as follows:

\begin{theorem}[Third Classification Theorem] \label{thm:DiscPoly}
The most general reducible discrete Markov process in the sense of
Definition \ref{def:OP}, has infinitesimal generator $\L$ given by
(\ref{eq:DiscL}). Its transitional probability density can be
expressed as
\begin{equation}\label{eq:DiscPolyU}
p_{T-t}(x, y) = \sum_{n=0}^N \frac{e^{\lambda_n(T-t)}}{d_n^2} Q_n(x)
Q_n(y) w(y).
\end{equation}
$\phi(x)$ must be of the form
\begin{equation}\label{eq:DiscPolyphi}
\phi(x) = B(x)+D(x)-\bar B(x)-\bar D(x),
\end{equation}
where $\bar B(x)$ and $\bar D(x)$ are the same rational functions as
$B(x)$ and $D(x)$ up to a multiplicative constant, but for different
parameters, and satisfy the lattice condition
\begin{equation}\label{eq:ConditionBD}
\bar B(x-1)\bar D(x) = B(x-1)D(x).
\end{equation}
The Laplace transform is given by the convergent series
(\mbox{$\displaystyle\prod_{k=1}^0=1$} by convention):
\begin{equation}\label{eq:DiscPolyG}
L_{T-t}(x, 1) =\prod_{k=1}^{x}\frac{D(k)}{\bar D(k)}\ \sum_{n=0}^N
e^{\bar\lambda_n(T-t)} z_n \bar Q_n(x).
\end{equation}
The coefficients $z_n$ are as follows:
\begin{equation}\label{eq:DiscPolyzn}
z_n = \frac{1}{\bar d_n^2}\sum_{x\in\Lambda_N}
\prod_{k=1}^x\frac{\bar D(k)}{D(k)}\ q(x)\bar Q_n(x) \bar w(x).
\end{equation}
\end{theorem}
\vspace{0.5cm} In our notation, $\{\bar Q_n(x)\}$ is the same set of
orthogonal polynomials as $\{Q_n(x)\}$ expect for the value of the
parameters. Section \ref{sec:DiscPoly} gives a proof of the latter
result, as well as explicit representations for processes based on
the Meixner, the dual Hahn and the Racah polynomials.\\

This paper is organized as follows: Section \ref{sec:Diffusion}
presents the proof of the first classification scheme in the
diffusion case and goes in more detail to provide a spectral
representation formula for transitional probability densities and
the Laplace transform. Section \ref{sec:Examples} contains a
discussion of the classical examples showing how the geometric
Brownian motion, the Ornstein-Uhlenbeck and the CIR process fit in
this classification scheme. Section \ref{sec:ContPoly} contains the
proof of an alternative classification scheme based on orthogonal
polynomials of the continuous series, along with interesting
examples. This discussion sets the premise for Section
\ref{sec:DiscPoly} where we extend the result to orthogonal
polynomials in the discrete series. Finally, in Section
\ref{sec:Limits} we discuss limiting relation and establish
connections between models corresponding to the discrete and the
continuous series, which is a useful result to construct numerical
very stable discretization schemes.

\newpage
\section{Classification Theorem for Diffusion Processes}
\label{sec:Diffusion}

In this section, we prove our first classification result, Theorem
\ref{thm:Diffusions}, in the diffusion case. We start by reviewing
some background notions concerning Fuchsian differential equations
and the so-called Bose invariants, and then proceed to the proof of
the theorem.

\subsection{Fuchsian Differential Equations}
Consider the second order partial differential equations for the
holomorphic function $F(z)$
\begin{equation}\label{eq:Fuchs}
F''(z) + p(z) F'(z) + q(z) F(z) = 0
\end{equation}
for some holomorphic functions $p(z), q(z)$.
\begin{definition}
Let $\alpha\in\C$ be an isolated singularity for the holomorphic
function $F(z)$. The singularity in $\alpha$ is called {\it regular}
if there is an exponent $\rho\in\C$ for which the function $(z -
\alpha)^{-\rho} F(z)$ admits a Laurent expansion with finitely many
negative powers around $z=\alpha$, i.e.
\begin{equation}
    F(z) = (z-\alpha)^\rho \sum^0_{n = -m} (z-\alpha)^n
\end{equation}
for some $m\in\N$. The point $\infty$ is a regular singularity of
the function $F(z)$ if $z=0$ is a regular singularity of the
function $F\left(\frac{1}{z}\right)$.
\end{definition}
In \cite{Gr1984}, Fuchs gives conditions on the coefficients $p(z)$
and $q(z)$ which ensure that solutions have only {\it regular}
singularities.
\begin{theorem}[Fuchs]
Let $F(z)$ be a solution of equation (\ref{eq:Fuchs}) with
singularities in the points $\alpha_1, \dots, \alpha_n$ and
$\infty$. Then these singularities are all regular if and only if
the functions $p(z)$ and $q(z)$ have the form
\begin{equation}
p(z) = {p_0(z)\over (z-\alpha_1)\dots (z-\alpha_n)}
\end{equation}
and
\begin{equation}
q(z) = {q_0(z)\over (z-\alpha_1)^2\dots (z-\alpha_n)^2}
\end{equation}
where $p_0(z)$ is a polynomial of order $(n-1)$ and $q_0(z)$ is a
polynomial of order $2 n - 2$.
\end{theorem}

An alternative expression for the coefficient $p(z)$ of an equation
with only regular singularities is
\begin{equation}
p(z) = \sum_{i = 1}^n {\delta_i\over z-\alpha_i}
\end{equation}
where the $\delta_i$, $i=1,..n$, are constants. In particular, we
have that
\begin{equation}
\exp\({1\over2} \int^z p(w) dw\) \ = \ C \prod_{i = 1}^n (
z-\alpha_i )^{\delta_i\over2}
\end{equation}
where $C$ is a constant. The function
\begin{equation}
\bar F(y) = \prod_{i = 1}^n ( y - \alpha_i )^{\delta_i\over2} F(y)
\end{equation}
solves the equation
\begin{equation}
\bar F''(y) + I(y) \bar F(y) = 0 \label{eq0}
\end{equation}
where
\begin{equation}
I(y) = -{1\over2}p'(y) + {1\over4} p(y)^2 + q(y).
\end{equation}

\begin{definition}
The function $I(y)$ is called the {\rm Bose invariant} of the
equation (\ref{eq:Fuchs}). Notice that $I(y)$ has the form
\begin{equation}
I(y) = {I_0(y)\over (y - \alpha_1)^2\dots (y - \alpha_n)^2}
\end{equation}
where $I_0(y)$ is a polynomial of order $2 n - 2$, without
restrictions on the coefficients.
\end{definition}

Let us focus again on the case $n=2$, assume that coefficients are
real and that only {\it real} linear fractional transformations are
allowed to move the regular singularities. In this situation we have
to distinguish between two different cases for Bose invariants:
\begin{itemize}
\item{} Case I ($\alpha_1 = 0, \; \alpha_2 = 1$):
\begin{equation}
I(y) = \frac{s_0 (1-y) + s_1 y + s_2 y(1-y)}{ y^2 (1-y)^2}
\end{equation}
\item{} Case II ($\alpha_1 = i, \; \alpha_2 = -i$):
\begin{equation}
I(y) = \frac{s_0 + s_1 y + s_2 y^2}{(y^2 + 1)^2}
\end{equation}
\end{itemize}
These cases reduce to the Gaussian hypergeometric equation
(\ref{eq:HypDiffEq}) for the function $_2F_1$ as is shown below.
Furthermore, special cases for the Bose invariant occur in the limit
when either $\alpha_1$ or $\alpha_2$ or both roots tend to $\infty$,
i.e.
\begin{itemize}
\item{} Case III ($\alpha_1 = 0, \ \alpha_2 = \infty$):
\begin{equation}
I(y) = \frac{s_0 + s_1 y + s_2 y^2}{y^2}
\end{equation}
\item{} Case IV ($\alpha_1 = \infty, \ \alpha_2 = \infty$):
\begin{equation}
I(y) = s_0 + s_1 y + s_2 y^2
\end{equation}
\end{itemize}
Case III reduces to the confluent hypergeometric equation $_1F_1$
and Case IV corresponds to the case of triple confluence at
infinity. Notice that the above four cases can all be captured by a
single expression as stated in the following:

\begin{proposition}\label{h_prop:BoseInv}
The Bose invariants corresponding to the Gaussian hypergeometric
function $_2F_1$ and to its confluent limit can be reduced to the
following normal form by means of a real valued linear fractional
transformation:
\begin{equation}
I(y) = \frac{Q(y)}{A(y)^2}
\end{equation}
where $A(y) \in \{y(1-y), y^2+1, y, 1\}$ and $Q$ is a polynomial in
y with $\deg Q \le 2$.
\end{proposition}

The first two cases correspond to three regular singularities at
distinct points. In these cases, solutions can be expressed through
Gaussian hypergeometric functions $_2F_1$. Fractional linear
transformations of the form
\begin{equation}
z\mapsto{az+b\over cz+d}
\end{equation}
where $a,b,c,d\in\C$ and $ad-bc\neq 0$, are one-to-one maps of the
extended complex line $\C\cup\infty$ into itself and map regular
singularities into regular singularities. By applying a fractional
linear transformation, one can map the singularities $\alpha_1$ and
$\alpha_2$ to 0 and 1, respectively. Furthermore, we have
transformations of the form
\begin{equation}
    F(z) \mapsto (z-\alpha_1)^{\rho_1} (z-\alpha_2)^{\rho_2} F(z).
\end{equation}
The combination of these two transformations allows one to reduce
any Fuchsian differential equation with three regular singular
points to the form
\begin{equation}\label{eq:HypDiffEq}
z(1-z) F''(z) + (\gamma - (1+\alpha+\beta) z) F'(z) - \alpha\beta
F(z) = 0.
\end{equation}
The function $_2F_1(\alpha,\beta; \gamma; z)$ is an elementary
solution of this equation along with $_2F_1(\alpha,\beta;
1+\alpha+\beta-\gamma; 1-z)$.

Case III corresponds to the limit when a regular singularity merges
with the regular singularity at $\infty$ while the other one stays
at 0. This limit can be obtained starting from the equation
corresponding to two coinciding singularities at 0, i.e. $\alpha_1 =
\alpha_2 = 0$:
\begin{equation}
F''(z) + {c_1 + c_2 z\over z^2} F'(z) + {c_3 + c_4 z + c_5 z^2\over
z^4} F(z) = 0.
\end{equation}
By applying the coordinate transformation $z\mapsto \frac{1}{z}$ we
find that
\begin{equation}\label{eq:conf0}
z^2 F''(z) + ((2 + c_2) z + c_1 z^2) F'(z) + (c_3 + c_4 z + c_5 z^2)
F(z) = 0
\end{equation}
By rescaling the independent variable $z$ and rescaling the function
so that $F(z) \mapsto e^{\rho z} f(\omega z)$ this equation reduces
to the Kummer differential equation
\begin{equation}\label{eq:HypDiffEq1F1}
z F''(z) + (\gamma - z) F'(z) - \alpha F(z) = 0
\end{equation}
which admits $_1F_1(\alpha;\gamma;z)$ as a solution. In alternative,
one can reduce equation (\ref{eq:conf0}) to the form
\begin{equation}
F''(z) + \( -\frac{1}{4} + \frac{\lambda}{z} + \frac{\frac{1}{4} -
\mu^2}{z^2} \) F(z) = 0
\end{equation}
which is called the {\it Whittaker differential equation}. The case
where all three singularities merge at $\infty$ is also interesting
and can be solved by rescaled confluent hypergeometric functions.

\subsection{Proof of the First Classification Theorem
\ref{thm:Diffusions}}

We start by presenting obvious facts about the transformations $T_i$
in Definition \ref{def:3op}.

\begin{proposition}\label{prop:inverse}
The transformations $T_i$ are invertible with respective inverse:
\begin{equation}
T_Z^{-1} = T_X\quad,\quad T_h^{-1} = T_\frac{1}{h}\quad,\quad
T_{\gamma^2}^{-1} = T_\frac{1}{\gamma^2}
\end{equation}
where $X:D_z\to D_x$ is the inverse of $Z(x)$.
\end{proposition}

\begin{proposition}\label{prop:commute}
The transformations $T_i$ commute with one another.
\end{proposition}
The proof of the First Classification Theorem \ref{thm:Diffusions}
for diffusion processes follows. The process $X_t$ has infinitesimal
generator
\begin{equation}
\L =  \frac{A(x)^2}{R(x)} \partII{}{x} + 2 \frac{h'(x)}{h(x)}\
\frac{A(x)^2}{R(x)}\ \partI{}{x}.
\end{equation}
To solve the eigenvalue problem $\L f =\lambda f$, define the
following left-multiplication $T_{\gamma^2}$ and gauge
transformation $T_h$:
\begin{equation}
T_{\gamma^2}\L = \frac{A(x)^2}{R(x)}\L\quad {\rm and}\quad T_h\L =
\frac{1}{h}\L h,
\end{equation}
so that the eigenvalue equation then transforms into
\begin{equation}
T_{\gamma^2}^{-1}T_h^{-1}: \L f=\lambda f\mapsto  \(\partII{}{x} +
\frac{S(x)}{A(x)^2}\)f(x) = \lambda\frac{R(x)}{A(x)^2}f(x).
\end{equation}
As $R(x),S(x)$ are second order polynomials and
$A(x)\in\{1,x,x(1-x),1+x^2\}$, the solution $f(x)$ can be expressed
as a hypergeometric function with Bose invariant
\begin{equation}
I(x) = \frac{S(x)-\lambda R(x)}{A(x)^2}.
\end{equation}
The same operations can be applied to the eigenvalue problem for the
Feynman-Kac operator, $\G \bar f = \bar\lambda\bar f$ with
\begin{equation}
\G = \L - \vartheta\phi(x),
\end{equation}
which leads to
\begin{equation}
T_{\gamma^2}^{-1}T_h^{-1}: \G \bar f=\bar\lambda\bar f\mapsto
\(\partII{}{x} + \frac{S(x)}{A(x)^2} -\vartheta\phi(x)\)\bar f(x) =
\bar\lambda\frac{R(x)}{A(x)^2}\bar f(x).
\end{equation}
For $\bar f(x)$ to be expressed as a hypergeometric function, we
require
\begin{equation}
\phi(x) = \frac{Q(x,\vartheta)}{\vartheta A(x)^2}
\end{equation}
for $Q(x,\vartheta)$ a second order polynomial in $x$. The Bose
invariant for the equation is
\begin{equation}
\tilde I(x) = \frac{S(x)-Q(x,\vartheta)-\bar\lambda R(x)}{A(x)^2}.
\end{equation}
The converse follows from the facts that the transformations are
invertible and the Bose invariants are both expressed in the most
general form.

\begin{remark}
The first transformation ($T_Y$: change of variable) has not been
used in the proof. $X_t$ is therefore the most general reducible
diffusion process only up to diffeomorphism.
\end{remark}

\subsection{Spectral Resolutions}\label{sub:Spectral}

Theorem \ref{thm:Diffusions} describes all processes with explicitly
solvable transitional probability density and Laplace transform for
the integral of the process. We wish now to give a closed form
expression for both quantities. The following lemma allows one to
determine the nature of the spectrum of the operators $\L$ and $\G$.
The nature of the spectrum is indeed based on the shape of the
Schr\"odinger potential, coming out of the eigenvalue problem once
reduced to a Schr\"odinger equation. Theorem \ref{thm:Schrodinger}
shows how this transformation operates on the kernel of the
semigroup generated by the Schr\"odinger operators and gives a
general spectral resolution of the operators $\L$ and $\G$.

\begin{lemma}\label{s_lem:Potentials}
Let $T = T_g\ T_Z\ T_h^{-1}$ where the diffeomorphism $Z:D_x\to D_z$
is given by $\displaystyle Z'(x)=\frac{\sqrt{R(x)}}{A(x)}$ with
inverse $X$ and the gauge transformation $T_g$ by $\displaystyle
g(z) = \(\frac{A(X(z))^2}{R(X(z))}\)^{1/4}$. Then the operators $\L$
and $\G$ reduce to the following Schr\"odinger operators:
\begin{eqnarray}\label{eq_schreq}
T \L &=& \partII{}{z} -
U_1(z)\equiv-\mathbb{H}_1\nonumber\\
T\G &=& \partII{}{z} - U_2(z)\equiv-\mathbb{H}_2
\end{eqnarray}
where the potentials are given by
\begin{eqnarray}
U(z) &=& \(\frac{g'}{g}\)^2 - \(\frac{g'}{g}\)' -
\frac{S(X(z))}{R(X(z))}\nonumber\\
\tilde U(z) &=& \(\frac{g'}{g}\)^2 - \(\frac{g'}{g}\)' -
\frac{S(X(z))-Q(X(z),\vartheta)}{R(X(z))}
\end{eqnarray}
and $'$ denotes the derivative with respect to $z$.
\end{lemma}

\begin{proof}
\begin{eqnarray}
T_g\ T_Z\ T_h^{-1}\L &=& T_g\ T_Z \(\frac{A(x)^2}{R(x)}\partII{}{x}
+ \frac
{S(x)}{R(x)} \)\nonumber\\
&=& T_g \(\partII{}{z} + \frac{Z''}{(Z')^2}\partI{}{z} +\frac
{S(X(z))}{R(X(z))} \)\nonumber\\
&=& \partII{}{z} +\frac{g''}{g} -2 \(\frac{g'}{g}\)^2 + \frac
{S(X(z))}{R(X(z))} \nonumber
\end{eqnarray}
and
$$ T_g\ T_Z\ T_h^{-1}\G = T_g\ T_Z\ T_h^{-1}\L - \vartheta\phi(X(z)).
$$
\end{proof}
The two Schr\"odinger operators $\mathbb{H}_1$ and $\mathbb{H}_2$
defined in the previous lemma have a spectral resolution,
\begin{equation}\label{eq:Schrodinger}
\mathbb{H}_i \Phi_\rho(z) = \rho\ \Phi_\rho(z),
\end{equation}
given by a complete set of normalized eigenfunctions $\Phi_\rho(z)$
for $\rho=-\lambda$ or $\rho=-\bar \lambda$, $i=1,2$ respectively.
The spectrum is in general a combination of a pure-point spectrum
$\sigma_{pp}(\mathbb{H}_i)$ and an absolutely continuous spectrum
$\sigma_{ac}(\mathbb{H}_i)$. The kernel of the semigroup generated
by the respective Schr\"odinger operators has the general form
\begin{equation}
e^{-(T-t)\mathbb{H}_i}(z_0, z_1) =
\sum_{\rho\in\sigma_{pp}(\mathbb{H}_i)} e^{-(T-t)\rho}
\Phi_\rho(z_0)\Phi_\rho^*(z_1) +
\int_{\rho\in\sigma_{ac}(\mathbb{H}_i)} e^{-(T-t)\rho}
\Phi_\rho(z_0)\Phi_\rho^*(z_1) dk(\rho)
\end{equation}
with $dk(\rho) = \frac{d\rho}{2\sqrt{\rho-U^-_i}}\displaystyle$ and
$U^-_i$ the lowest limit of the potential $U_i(z)$ as $z$ tends to
the boundaries of the domain $D_z$.

\begin{theorem}\label{thm:Schrodinger}
The transitional probability density and the Laplace transform of
any reducible process described in Theorem \ref{thm:Diffusions} by
the operator $\L$ and $\G$ is related to the kernels of the
semigroups generated by the respective Schr\"odinger operators as
follows:
\begin{eqnarray}\label{eq:LGkernel}
e^{(T-t)\L}(x, y) & =& \frac{h(y)}{h(x)}\
\(\frac{A(x)^2}{R(x)}\)^{1/4} \(\frac{R(y)}{A(y)^2}\)^{3/4}
e^{-(T-t)\mathbb{H}_1}\big(Z(x), Z(y)\big),
\nonumber\\
e^{(T-t)\G}(x, y) & =& \frac{h(y)}{h(x)}\
\(\frac{A(x)^2}{R(x)}\)^{1/4} \(\frac{R(y)}{A(y)^2}\)^{3/4}
e^{-(T-t)\mathbb{H}_2}\big(Z(x), Z(y)\big).
\end{eqnarray}
\end{theorem}

\begin{proof}
The transformation $T$ defined in Lemma \ref{s_lem:Potentials}
extends to the kernels as follows:
\begin{eqnarray}
e^{(T-t)\L}(x, y) &=& e^{-(T-t)\ T_h T_Z^{-1} T_g^{-1}
\mathbb{H}_1}(x, y)\nonumber\\ &=& \frac{h(y)}{h(x)}\ \frac{dZ}{dy}\
\frac{g(Z(x))}{g(Z(y))}\ e^{-(T-t)\mathbb{H}_1}\big(Z(x),
Z(y)\big)\nonumber
\end{eqnarray}
and similarly for $\G$ and $\mathbb{H}_2$. Recall that
$\displaystyle g(Z(x)) =\(\frac{A(x)^2}{R(x)}\)^{1/4}$ and the
Jacobian of $T_Z$ is given by
$\displaystyle\frac{dZ}{dy}=\frac{\sqrt{R(y)}}{A(y)}$, which
concludes the proof.
\end{proof}

\newpage
\section{Examples of Solvable Diffusions}\label{sec:Examples}

In this section, we show that the transitional probability density
and the Laplace transform for the integral of the geometric Brownian
motion, the Ornstein-Uhlenbeck process and the CIR process arise as
corollaries of the First Classification Theorem
\ref{thm:Diffusions}.

\subsection{The geometric Brownian motion}\label{sub:GBM}

\begin{definition}
The geometric Brownian motion is defined by the solution of the
following stochastic differential equation:
\begin{equation}
dX_t = \mu X_t dt + \sigma X_t dW_t
\end{equation}
with initial condition $X_{t=0}=x_0$.
\end{definition}

\begin{corollary}
The transitional probability density for the geometric Brownian
motion is given by the following formula:
\begin{equation}
p_{T-t}(x, y)= \frac{1}{y\sqrt{2\pi\sigma^2(T-t)}}\
\exp\(-\frac{\(\ln(\frac{y}{x})-(\mu-\frac{\sigma^2}{2})(T-t)\)^2}
{2\sigma^2(T-t)}\).
\end{equation}
The Laplace transform is explicitly solvable if and only if
\begin{equation}
\phi(x) = \frac{\sigma^2}{2\vartheta}
\(\frac{\mu}{\sigma^2}\Big(1-\frac{\mu}{\sigma^2}\Big)-t_0 - t_1 x +
t_2 x^2\),
\end{equation}
where $t_0, t_1\in\R$ and $t_2>0$ could depend on $\vartheta$. It is
then given in terms of the Laguerre polynomials $L_n^{(\delta)}$ and
the Whittaker function $M_{\lambda,\mu}$ (by convention,
$\displaystyle\sum_{n=0}^N=0$ if $N<0$):
\begin{eqnarray}\label{eq:GBMG}
L_{T-t}(x, \vartheta) &=& x^{-\frac{\mu}{\sigma^2}+\frac{1}{2}}\
e^{-\sqrt{t_2}x}\ \sum_{n=0}^N z_n e^{(T-t)\lambda_n}
x^\frac{\delta_n}{2}
L_n^{(\delta_n)}(2\sqrt{t_2}x)\nonumber\\
&+& x^{-\frac{\mu}{\sigma^2}} \int_{0}^\infty z_k
e^{-(T-t)(k^2+U_-)}
M_{\frac{t_1}{2\sqrt{t_2}},\frac{\delta_k}{2}}(2\sqrt{t_2}x) dk
\end{eqnarray}
where $N\equiv\lceil\frac{t_1}{\sqrt{t_2}}-\frac{1}{2}\rceil$
($\lceil t\rceil$ denotes the integer part of $t$),
$\delta_n\equiv-2n-1+\frac{t_1}{\sqrt{t_2}}$,
$U_-\equiv\frac{\sigma^2}{2}(\frac{1}{4} - t_0)$ and $\delta_k
\equiv i\sqrt\frac{8}{\sigma^2}k$. For $n=0,1,\ldots,N$, the
discrete eigenvalues are given by
\begin{equation}\label{eq:GBMev}
\lambda_n = \frac{\sigma^2}{2} \(\(n+\frac{1}{2}
-\frac{t_1}{2\sqrt{t_2}}\)^2 +t_0-\frac{1}{4}\).
\end{equation}
The coefficients $z_n$ and $z_k$ are respectively
\begin{equation}
z_n = \frac{(2\sqrt{t_2})^{\delta_n}}{\Gamma(\delta_n)}
\int_0^\infty q(x)\
x^{\frac{\mu}{\sigma^2}+\frac{\delta_n}{2}-\frac{3}{2}}
e^{-\sqrt{t_2}x} L_n^{(\delta_n)}(2\sqrt{t_2}x)dx
\end{equation}
and
\begin{equation}
z_k = \frac{1}{2\pi}\sqrt\frac{2}{\sigma^2t_2} \int_0^\infty q(x)\
x^{\frac{\mu}{\sigma^2}-2} M_{\frac{t_
1}{2\sqrt{t_2}},\frac{-\delta_k}{2}}(2\sqrt{t_2}x) dx.
\end{equation}
\end{corollary}

\begin{proof}
The infinitesimal generator is
\begin{equation}
\L = \frac{\sigma^2}{2}x^2\partII{}{x} + \mu x\partI{}{x}
\end{equation}
whereas the Feynman-Kac operator has the form
\begin{equation}
\G = \frac{\sigma^2}{2}x^2\partII{}{x} + \mu x\partI{}{x} +
\vartheta\phi(x).
\end{equation}
This case fits the classification scheme in Theorem
\ref{thm:Diffusions} if one selects
\begin{equation}
A(x)=x\quad,\quad R(x)=\frac{2}{\sigma^2}\quad,\quad \frac{h'}{h} =
\frac{\mu}{\sigma^2x}.
\end{equation}
This choice sets the shape of the polynomial $S(x)$, from
(\ref{eq:T1v}) in Theorem \ref{thm:Diffusions}, to be
\begin{equation}
S(x) =
\frac{1}{\sigma^2}\(2\xi+\mu\big(1-\frac{\mu}{\sigma^2}\big)\),
\end{equation}
which in turn defines the Bose invariant as
\begin{equation}
I(x) = \frac{2(\xi-\lambda)
+\mu\big(1-\frac{\mu}{\sigma^2}\big)}{\sigma^2x^2}.
\end{equation}
The Schr\"odinger potential, given by
\begin{equation}
U(z) = \frac{\sigma^2}{8}-\xi
-\frac{\mu}{2}\(1-\frac{\mu}{\sigma^2}\),
\end{equation}
is constant for all $z$. Hence the spectrum is absolutely
continuous. The normalized eigenfunctions for the Schr\"odinger
equations (\ref{eq:Schrodinger}) are
\begin{equation}
\Phi_\rho(z) = \frac{1}{\sqrt{2\pi}}\ e^{\pm ik(\rho)z}
\end{equation}
where
$k^2=\rho-\frac{\sigma^2}{8}+\frac{\mu}{2}\(1-\frac{\mu}{\sigma^2}\)$.
Theorem \ref{thm:Schrodinger} yields the kernel for the semigroup
generated by the operator $\L$:
\begin{equation}
e^{(T-t)\L}(x,y)= \frac{1}{y\sqrt{2\pi\sigma^2(T-t)}}\
\exp\(-\frac{\(\ln(\frac{y}{x})-(\mu-\frac{\sigma^2}{2})(T-t)\)^2}
{2\sigma^2(T-t)}\) 
\end{equation}
which is $p_{T-t}(x, y)$ and in which one recognizes the
transitional probability density of the geometric Brownian motion.

The kernel of the semigroup generated by $\G$ is more general, since
there is no restrictions on the polynomial $Q(x, \vartheta)$, which,
for $t_0,t_1,t_2\in\R$, can be written in the form
\begin{equation}
Q(x, \vartheta) = S(x) - t_0-t_1x+t_2x^2.
\end{equation}
The Bose invariant in this case is
\begin{equation}
\tilde I(x) =
-t_2+\frac{t_1}{x}+\frac{t_0-\frac{2\lambda}{\sigma^2}}{x^2}
\end{equation}
which gives rise to two independent solutions to the eigenvalue
problem $u'' + Iu=0$:
\begin{eqnarray}
u_\pm(x) &=&
M_{\frac{t_1}{2\sqrt{t_2}},\pm\frac{\delta_0}{2}}\(2\sqrt{t_2}x\)
\nonumber\\
&=&(2\sqrt{t_2}x)^{\frac{\pm\delta_0+1}{2}} e^{-\sqrt{t_2}x}
{}_1F_1\Bigg(\begin{matrix}
\frac{\pm\delta_0+1}{2}-\frac{t_1}{2\sqrt{t_2}}\\
\pm\delta_0+1\end{matrix}\ \Bigg\vert\ 2\sqrt{t_2}x\Bigg)
\end{eqnarray}
where $\frac{\delta_0^2}{4} = \frac{1}{4}-t_0
+\frac{2\lambda}{\sigma^2}$. The Schr\"odinger potential, given by
\begin{equation}
\tilde U(z) = \frac{\sigma^2}{2}\(\frac{1}{4} - t_0
-t_1e^{\sqrt{\frac{\sigma^2}{2}}z}
+t_2e^{2\sqrt{\frac{\sigma^2}{2}}z}\),
\end{equation}
has no singularities and is bounded from below if $t_2>0$. It
indicates that the spectrum is not strictly discrete, since $\tilde
U(z)\to\infty$ as $z\to\infty$ but $\tilde U(z)\to
U_-\equiv\frac{\sigma^2}{2}(\frac{1}{4} - t_0)$ as $z\to-\infty$.
Hence, the description of the spectrum separates in two cases:
\begin{itemize}
\item If $t_1\le0$, then $\tilde U(z)$ is monotonously increasing
on $D_z=\R$ and the spectrum is continuous: $\rho=-\lambda\ge U_-$.

\item If $t_1>0$, then $\tilde U(z)$ has a minimum, $U_0\equiv U_-
-\frac{\sigma^2}{2}\frac{t_1^2}{4t_2}$, at
$z=\sqrt\frac{2}{\sigma^2}\ln\frac{t_1}{2t_2}$ and the spectrum is
discrete for $U_0<\rho<U_-$ and continuous for $\rho\ge U_-$.

\end{itemize}
Both cases can however be solved simultaneously. The solution to the
Schr\"odinger equation, $\Phi_\rho(z)$, is related to $u_\pm(x)$,
via the diffeomorphism $x=X(z)=e^{\sqrt\frac{\sigma^2}{2}z}$ and the
gauge transformation $g(z)=(\frac{\sigma^2}{2})^\frac{1}{4}
e^{\sqrt\frac{\sigma^2}{8}z}$, as follows:
\begin{equation}
\Phi_\rho(z) = g^{-1}(z) u_\pm(X(z)).
\end{equation}
The asymptotic behavior as $z\to-\infty$, $\Phi_\rho(z)\approx
\frac{1}{\sqrt{\pi}}e^{\pm ikz}$ with $k=\sqrt{\rho-U_-}\ge0$,
enforces the normalization condition for $\Phi_\rho(z)$ to the
following:
\begin{equation}
\Phi_\rho(z) = \frac{1}{\sqrt{\pi}}
(2\sqrt{t_2})^{-\frac{1\pm\delta_0}{2}}
e^{-\sqrt\frac{\sigma^2}{8}z} u_\pm(X(z)).
\end{equation}
The continuous spectrum appears for $\rho$ greater than the lowest
limit $U_-$, implying that $\delta_k\equiv\delta_0(k)=
i\sqrt\frac{8}{\sigma^2}k$ is imaginary. Whether or not the discrete
part of the spectrum has an infinite number of discrete levels
depends on the asymptotic behavior of the potential. As
$z\to-\infty$, $x\to0$ and $\tilde U(Z(x))$ develops up to second
order:
\begin{equation}
\tilde U(Z(x)) = \frac{\sigma^2}{2}\(\frac{1}{4}-t_0-t_1 x+t_2
x^2\).
\end{equation}
Since $t_1>0$ for the discrete spectrum, the convergence to the
limit $U_-$ as $x\to0$ is faster than $x^2$, which implies that
there is only a finite number of bound states.

The hypergeometric functions in the solutions $u_\pm(x)$ reduce to
polynomials for $x\in[0,\infty)$ if respectively
\begin{equation}
\frac{\pm\delta_0(n)+1}{2}-\frac{t_1}{2\sqrt{t_2}} = -n\in-\N.
\end{equation}
We set $\delta_n\equiv\delta_0(n)$ to emphasize the dependance on
$n$. The solution $u_-(x)$ is not $L^2$-normalizable, whereas
$u_+(x)$ is given in terms of the Laguerre polynomials
\begin{equation}
u_n(x) = C (2\sqrt{t_2}x)^{\frac{\delta_n+1}{2}} e^{-\sqrt{t_2}x}
L_n^{(\delta_n)}(2\sqrt{t_2}x)
\end{equation}
with normalization constant $C$. The normalization condition
$\displaystyle\int_{-\infty}^\infty\vert\Phi_n(z)\vert^2dz=1$ fixes
the constant $\displaystyle C=\sqrt\frac{\sigma^2}
{4\sqrt{t_2}\Gamma(\delta_n)}$ (cf. \cite{PBM1986}, p. 462). The
discrete eigenvalues $\rho=-\lambda_n$ are given by
(\ref{eq:GBMev}), where $n$ is restricted to $\{0,1,\ldots,N\}$ by
the condition $-\lambda_n>U_0$, and
$N\equiv\lceil\frac{t_1}{\sqrt{t_2}}-\frac{1}{2}\rceil$.

Hence, from Theorem \ref{thm:Schrodinger}, the kernel of the
semigroup generated by $\G$ is given by the following spectral
resolution:
\begin{eqnarray}
e^{(T-t)\G}(x,y)& = &x^{-\frac{\mu}{\sigma^2}+\frac{1}{2}}
y^{\frac{\mu}{\sigma^2}-\frac{3}{2}}\ e^{-\sqrt{t_2}(x+y)}\\
&&\cdot\sum_{n=0}^N e^{(T-t)\lambda_n}
\frac{(2\sqrt{t_2})^{\delta_n}}{\Gamma(\delta_n)}
(xy)^\frac{\delta_n}{2}
L_n^{(\delta_n)}(2\sqrt{t_2}x) L_n^{(\delta_n)}(2\sqrt{t_2}y)\nonumber\\
&+&\frac{1}{2\pi}\sqrt\frac{2}{\sigma^2t_2}
x^{-\frac{\mu}{\sigma^2}}
y^{\frac{\mu}{\sigma^2}-2} \nonumber\\
&&\cdot\int_{0}^\infty
e^{-(T-t)(k^2+U_-)}
M_{\frac{t_1}{2\sqrt{t_2}},\frac{\delta_k}{2}}(2\sqrt{t_2}x)
M_{\frac{t_1}{2\sqrt{t_2}},\frac{-\delta_k}{2}}(2\sqrt{t_2}y) dk.
\nonumber
\end{eqnarray}
Finally the Laplace transform of the latter kernel is given by
integration over $D_x=[0,\infty)$ and yields (\ref{eq:GBMG}).
\end{proof}

\subsection{The Ornstein-Uhlenbeck process}\label{sub:OU}
In this subsection, we restrict the framework to the affine models,
i.e. we set
\begin{equation}
q(x) = \exp\big(\omega \phi(x)\big)
\end{equation}
for some $\omega\in\R$. We show that in the special case of the
Ornstein-Uhlenbeck process, both the transitional probability
density and the Laplace transform can be expressed as summations
over Hermite polynomials.

\begin{definition}
The Ornstein-Uhlenbeck process is defined by the solution of the
following stochastic differential equation:
\begin{equation}
dX_t = (a-bX_t) dt + \sigma dW_t
\end{equation}
with $b>0$ and initial condition $X_{t=0}=x_0$.
\end{definition}

\begin{corollary}
The transitional probability density is given by the following
formula:
\begin{equation}\label{eq:OUU}
p_{T-t}(x, y) =
\sqrt{\frac{b}{\sigma^2\pi}}\(1-e^{-2b(T-t)}\)^{-1/2}
\exp\[-\frac{\big(z(y)-z(x)e^{-b(T-t)}\big)^2}{1-e^{-2b(T-t)}}\]
\end{equation}
where $z(x) = \sqrt{\frac{b}{\sigma^2}}\(x-\frac{a}{b}\)$. The
Laplace transform is explicitly solvable if and only if
\begin{equation}
\phi(x) = \frac{\sigma^2}{2\vartheta}
\(\frac{b}{\sigma^2}-\frac{a^2}{\sigma^4}-t_0
+2\frac{ab-\a\b}{\sigma^4} x + \frac{\b^2-b^2}{\sigma^4} x^2\),
\end{equation}
with $t_0, \a,\b\in\R$ and could depend on $\vartheta$. It is of the
quadratic form
\begin{equation}\label{eq:OUG}
L_{T-t}(x, \vartheta,\omega) = e^{m(T-t) - n(T-t)x - l(T-t)x^2}
\end{equation}
where the functions $m(\tau), n(\tau)$ and $l(\tau)$ are as follows
\begin{eqnarray}
m(\tau) &=& \frac{1}{2}\ln\(\frac{\b}{\sigma^2}\)
-\frac{1}{2}\ln\(p-(p-\frac{\b}{\sigma^2})e^{-2\b\tau}\)
-\frac{\a^2}{\sigma^2\b} -
\frac{\tau}{2}\(\b-\frac{\a^2}{\sigma^2}-\sigma^2t_0\)\nonumber\\
&&+\frac{\omega}{2\vartheta}\(b-\frac{a^2}{\sigma^2}-\sigma^2t_0\)
+p\(q+\frac{\a}{\b}\)^2 - \frac{2pq\frac{\a}{\sigma^2}e^{-\b\tau}
+\((p-\frac{\b}{\sigma^2})\frac{\a^2}{\b\sigma^2}
+pq^2\frac{\b}{\sigma^2}\)e^{-2\b\tau}}
{p-(p-\frac{\b}{\sigma^2})e^{-2\b\tau}}
\nonumber\\
n(\tau)&=& -\frac{\a-a}{\sigma^2} - \frac{2pq\frac{\b}{\sigma^2}
e^{-\b\tau} +
(p-\frac{\b}{\sigma^2})2\frac{\a}{\sigma^2}e^{-2\b\tau}}
{p-(p-\frac{\b}{\sigma^2})e^{-2\b\tau}} \nonumber\\
l(\tau)&=& \frac{\b-b}{2\sigma^2} +
\frac{(p-\frac{\b}{\sigma^2})\frac{\b}{\sigma^2} e^{-2\b\tau}}
{p-(p-\frac{\b}{\sigma^2})e^{-2\b\tau}}\nonumber
\end{eqnarray}
with
\begin{equation}
p= \frac{(\b+b)(\omega(b-\b)+\vartheta)}{2\sigma^2\vartheta}\nonumber\\
\quad,\quad q= \frac{\vartheta(a+\a)+\omega(ab-\a\b)}
{(\b+b)(\vartheta-\omega(\b-b))}-\frac{\a}{\b}.
\end{equation}
\end{corollary}

\begin{proof}
The infinitesimal generator of the Ornstein-Uhlenbeck process is
\begin{equation}
\L = \frac{\sigma^2}{2}\partII{}{x} + (a-bx)\partI{}{x}
\end{equation}
and the Feynman-Kac operator has the form
\begin{equation}
\G = \frac{\sigma^2}{2}\partII{}{x} +(a-bx)\partI{}{x} -
\vartheta\phi(x).
\end{equation}
This case fits the classification scheme in Theorem
\ref{thm:Diffusions} if one selects
\begin{equation}
A(x)=1\quad,\quad R(x)=\frac{2}{\sigma^2}\quad,\quad h(x) =
e^{\frac{a}{\sigma^2}x-\frac{b}{2\sigma^2}x^2}.
\end{equation}
This choice implies that the shape of the polynomial $S(x)$, from
(\ref{eq:T1v}) in Theorem \ref{thm:Diffusions}, must be
\begin{equation}
S(x) = \frac{1}{\sigma^2}(2\xi+b-\frac{a^2}{\sigma^2}) +
\frac{2ab}{\sigma^4}x - \frac{b^2}{\sigma^4}x^2,
\end{equation}
which in turn defines the Bose invariant as
\begin{equation}
I(x) = \frac{1}{\sigma^2}(2\xi-2\lambda+b-\frac{a^2}{\sigma^2}) +
\frac{2ab}{\sigma^4}x - \frac{b^2}{\sigma^4}x^2.
\end{equation}
The Schr\"odinger potential, given by
\begin{equation}
U(z) = -\xi -\frac{b}{2}
+\frac{a^2}{2\sigma^2}-\frac{ab}{\sqrt{2\sigma^2}}z
+\frac{b^2}{4}z^2
\end{equation}
with $z=\sqrt{\frac{2}{\sigma^2}}x$, goes to $\infty$ as
$z\to\pm\infty$. Hence, the spectrum displays discrete eigenvalues
of the form $\lambda_n=-bn$ for $n\in\N$. The corresponding
normalized eigenfunctions that satisfy the Schr\"odinger equation
(\ref{eq:Schrodinger}) are
\begin{equation}
\Phi_n(Z(x)) =\frac{1}{\sqrt{n!2^n}}\(\frac{b}{2\pi}\)^{1/4}
e^{-\frac{b}{2\sigma^2}\(x-\frac{a}{b}\)^2}
H_n\(\sqrt{\frac{b}{\sigma^2}}\big(x-\frac{a}{b}\big)\).
\end{equation}
Theorem \ref{thm:Schrodinger} yields the kernel for the semigroup
generated by $\L$ as a summation over Hermite polynomials. Using
Mehler's formula (cf. \cite{PBM1986}, p. 710) and the notation $z(x)
= \sqrt{\frac{b}{\sigma^2}}\(x-\frac{a}{b}\)$, the kernel re-sums
into
\begin{equation}
e^{(T-t)\L}(x,y)=
\sqrt{\frac{b}{\sigma^2\pi}}\(1-e^{-2b(T-t)}\)^{-1/2}
\exp\[-\frac{\big(z(y)-z(x)e^{-b(T-t)}\big)^2}{1-e^{-2b(T-t)}}\],
\end{equation}
which is $p_{T-t}(x, y)$ and in which one recognizes the probability
density of the Ornstein-Uhlenbeck process.

For convenience and without loss of generality, we set the form of
the polynomial $Q(x,\vartheta)$ to
\begin{equation}
Q(x,\vartheta) = S(x) -  t_0 -
\frac{2\a\b}{\sigma^4}x+\frac{\b^2}{\sigma^4}x^2.
\end{equation}
where $t_0,\a,\b$ could depend on $\vartheta$. The Bose invariant in
this case is
\begin{equation}
\tilde I(x) = t_0-\frac{2\lambda}{\sigma^2}+
\frac{2\a\b}{\sigma^4}x-\frac{\b^2}{\sigma^4}x^2.
\end{equation}
The Schr\"odinger potential, given by
\begin{equation}
\tilde U(z) = -\frac{t_0}{2} - \frac{\a\b}{\sigma^2}
\sqrt\frac{\sigma^2}{2}z +\frac{\b^2}{4}z^2
\end{equation}
is very similar to $U(z)$ and indicates that the spectrum is again
discrete since $\tilde U(z)\to\infty$ as $z\to\pm\infty$. The
solution $\Phi_n(y)$ is given in terms of the Hermite polynomials
\begin{equation}
\Phi_n(Z(x)) =\frac{1}{\sqrt{n!2^n}}\(\frac{\b}{2\pi}\)^{1/4}
e^{-\frac{\b}{2\sigma^2}\(x-\frac{\a}{\b}\)^2}
H_n\(\sqrt{\frac{\b}{\sigma^2}}\big(x-\frac{\a}{\b}\big)\).
\end{equation}
for the eigenvalues $\lambda_n = -\b n -
\frac{\b}{2}+\frac{\a^2}{2\sigma^2} +\frac{\sigma^2}{2}t_0$. From
Theorem \ref{thm:Schrodinger}, the kernel of the semigroup generated
by $\G$ is given by the following spectral resolution:
\begin{eqnarray}
e^{(T-t)\G}(x,y) &=& \sqrt\frac{\b}{\sigma^2\pi}
\exp\(-\frac{\a^2}{\sigma^2\b}+\frac{y}{\sigma^2}(\a+a) -
\frac{y^2}{2\sigma^2}(\b+b)+\frac{x}{\sigma^2}(\a-a) -
\frac{x^2}{2\sigma^2}(\b-b)\)\nonumber\\
&&\cdot\sum_{n=0}^\infty \frac{e^{\lambda_n(T-t)}}{2^nn!}
H_n\(\sqrt\frac{\b}{\sigma^2}\big(x-\frac{\a}{\b}\big)\)
H_n\(\sqrt\frac{\b}{\sigma^2}\big(y-\frac{\a}{\b}\big)\).
\end{eqnarray}
The integration of the latter yields the Laplace transform for $q(x)
= \exp(\omega x)$ as a convergent series in terms of the Hermite
polynomials which re-sums to the formula (\ref{eq:OUG}) (cf.
\cite{PBM1986}, p. 488(16) and p. 710(1)).
\end{proof}

\begin{remark}
Computing the Laplace transform $L$ in the case where $\phi(x)$ is
affine, i.e. $\phi(x)=x$, is a direct consequence of the previous
corollary. Setting the parameters to
\begin{equation}
t_0=\frac{b}{\sigma^2} - \(\frac{a}{\sigma^2}\)^2\quad,\quad
\a=a-\frac{\sigma^2\vartheta}{b}\quad,\quad \b=b
\end{equation}
proves the following proposition:
\begin{proposition}
The Laplace transform for the affine Ornstein-Uhlenbeck process is
as follows:
\begin{equation}\label{eq:OUGaffine}
L_{T-t}(x,\vartheta,\omega) = e^{m(T-t) - n(T-t)x}
\end{equation}
where
\begin{eqnarray}
n(\tau) &=& \frac{\vartheta-(\vartheta+\omega b)e^{-b\tau}}{b}\nonumber\\
m(\tau) &=&
\frac{\big(n(\tau)+\omega-\vartheta\tau\big)(ab-\vartheta\frac{\sigma^2}{2})}{b^2}
- \frac{\sigma^2}{4b}\big(n(\tau)^2-\omega^2\big).
\end{eqnarray}
\end{proposition}
\end{remark}

\subsection{The CIR process}\label{sub:CIR}
In this subsection, we focuss again on the affine models, i.e. we
set
\begin{equation}
q(x) = \exp(\omega x)
\end{equation}
for some $\omega\in\R$. We show how to derive from the First
Classification Theorem \ref{thm:Diffusions}, the transitional
probability density for the CIR process and the Laplace transform in
the affine case. The main tool is the use of the Laguerre
polynomials as eigenfunctions for both the infinitesimal generator
and the Feynman-Kac operator.

\begin{definition}
The CIR process is defined as the solution of the following
stochastic differential equation on $D_x=\R_+$:
\begin{equation}
dX_t = (a-bX_t) dt + \sigma\sqrt{X_t}dW_t
\end{equation}
with $a,b>0$ and initial condition $X_{t=0}=x_0$.
\end{definition}

\begin{corollary}
The transitional probability density of the CIR process is given in
terms of the modified Bessel function $I_\alpha$ as follows:
\begin{eqnarray}\label{eq:CIRU}
p_{T-t}(x, y)&=& c\(\frac{y
e^{b(T-t)}}{x}\)^{\frac{1}{2}\(\frac{2a}{\sigma^2}-1\)}
\exp\[-c\(y+x e^{-b(T-t)}\)\] \nonumber
I_{\frac{2a}{\sigma^2}-1}\(2c\sqrt{x y e^{-b(T-t)}}\)\\
\end{eqnarray}
with $c\equiv c(T-t) = \frac{2b}{\sigma^2}(1-e^{-b(T-t)})^{-1}$. The
Laplace transform is computable in closed form if and only if
\begin{equation}
\phi(x) = \(\frac{a}{2}\big(1-\frac{a}{\sigma^2}\big)
-\frac{\sigma^2t_0}{2}\)\frac{1}{x}
+\(\frac{ab}{\sigma^2}-\frac{\sigma^2t_1}{2}\) + \(
\frac{\sigma^2t_2}{2}-\frac{b^2}{2\sigma^2}\)x,
\end{equation}
where $t_0, t_1,t_2\in\R$ could depend on $\vartheta$.

In the affine case, where $\phi(x)=x$, the Laplace transform is of
the closed form
\begin{equation}\label{eq:CIRG}
L_{T-t}(x,\vartheta,\omega) = e^{m(T-t) - n(T-t)x}
\end{equation}
where
\begin{eqnarray}
m(\tau) &=& \frac{2a}{\sigma^2}\ln\[\frac{\b
e^{b\tau/2}}{\b\cosh(\frac{\b\tau}{2}) +
(b-\omega\sigma^2)\sinh(\frac{\b\tau}{2})}\]\nonumber\\
n(\tau) &=& -\omega+\frac{\b^2-(b-\omega\sigma^2)^2}{\sigma^2}\
\frac{\sinh(\frac{\b\tau}{2})}{\b\cosh(\frac{\b\tau}{2}) +
(b-\omega\sigma^2)\sinh(\frac{\b\tau}{2})}
\end{eqnarray}
and $\b=\sqrt{2\vartheta\sigma^2+b^2}$.
\end{corollary}

\begin{proof}
The infinitesimal generator of the CIR process is
\begin{equation}
\L = \frac{\sigma^2}{2}x\partII{}{x} + (a-bx)\partI{}{x}
\end{equation}
and the Feynman-Kac operator has the form
\begin{equation}
\G = \frac{\sigma^2}{2}x\partII{}{x} +(a-bx)\partI{}{x} -
\vartheta\phi(x).
\end{equation}
This case fits the classification scheme in Theorem
\ref{thm:Diffusions} if one selects
\begin{equation}
A(x)=x\quad,\quad R(x)=\frac{2x}{\sigma^2}\quad,\quad h(x) =
x^{a/\sigma^2}e^{-\frac{b}{\sigma^2}x}.
\end{equation}
This choice sets the form of the polynomial $S(x)$, from
(\ref{eq:T1v}) in Theorem \ref{thm:Diffusions}, to be
\begin{equation}
S(x) = \frac{a}{\sigma^2}\(1-\frac{a}{\sigma^2}\)
+\frac{2}{\sigma^2}\(\xi+\frac{ab}{\sigma^2}\)x -
\frac{b^2}{\sigma^4}x^2,
\end{equation}
which in turn defines the Bose invariant as
\begin{equation}
I(x) = - \frac{b^2}{\sigma^4}
+\frac{2}{\sigma^2}\(\xi-\lambda+\frac{ab}{\sigma^2}\)\frac{1}{x}
+\frac{a}{\sigma^2}\(1-\frac{a}{\sigma^2}\)\frac{1}{x^2}.
\end{equation}
The Schr\"odinger potential, given for $z=\sqrt\frac{8x}{\sigma^2}$
by
\begin{equation}
U(z) = \(\frac{3}{4}
-\frac{4a}{\sigma^2}\(1-\frac{a}{\sigma^2}\)\)\frac{1}{z^2}
-\xi-\frac{ab}{\sigma^2}+ \frac{b^2}{16} z^2
\end{equation}
goes to $\infty$ as $z\to\infty$ and has a singularity at $z=0$. The
asymptotic behavior as $z\to0$ implies that the Schr\"odinger
operator has a self-adjoint closure with completely discrete
spectrum. We have assumed
$\frac{4a}{\sigma^2}\(1-\frac{a}{\sigma^2}\)\le1$ which is always
satisfied. Hence, the spectrum presents discrete eigenvalues of the
form $\lambda_n=-bn$ for $n\in\N$. The corresponding normalized
eigenfunctions that satisfy the Schr\"odinger equation
(\ref{eq:Schrodinger}) are given in terms of the Laguerre
polynomials
\begin{equation}
\Phi_n(Z(x)) =\sqrt\frac{n!}{\Gamma(n+\frac{2a}{\sigma^2})}
\(\frac{\sigma^2}{2x}\)^{1/4} \(\frac{2bx}{\sigma^2}\)^{a/\sigma^2}
e^{-\frac{b}{\sigma^2}x}
L_n^{(\frac{2a}{\sigma^2}-1)}\(\frac{2b}{\sigma^2}x\).
\end{equation}
Theorem \ref{thm:Schrodinger} yields the kernel for the semigroup
generated by the operator $\L$ as a summation over Laguerre
polynomials, which can be re-summed into (cf. \cite{PBM1986}, p.
705(7)):
\begin{eqnarray}
e^{(T-t)\L}(x,y)&=& c\(\frac{y
e^{b(T-t)}}{x}\)^{\frac{1}{2}\(\frac{2a}{\sigma^2}-1\)}
\exp\[-c\(y+x e^{-b(T-t)}\)\] \nonumber
I_{\frac{2a}{\sigma^2}-1}\(2c\sqrt{x y e^{-b(T-t)}}\)\\
\end{eqnarray}
with $c\equiv c(T-t) = \frac{2b}{\sigma^2}(1-e^{-b(T-t)})^{-1}$ and
in which one recognizes the transitional probability density of the
CIR process.

For the Laplace transform $L$, we specialize to the case where
$\phi(x)$ is affine, i.e. $\phi(x)=x$. Therefore, (\ref{eq:phi})
sets the shape of the polynomial $Q(x,\vartheta)$ to
\begin{equation}
Q(x,\vartheta) = S(x) +\frac{2\vartheta}{\sigma^2}x^2.
\end{equation}
The Bose invariant in this case is
\begin{equation}
\tilde I(x) =  - \frac{b^2}{\sigma^4}-\frac{2\vartheta}{\sigma^2}
+\frac{2}{\sigma^2}\(\xi-\lambda+\frac{ab}{\sigma^2}\)\frac{1}{x}
+\frac{a}{\sigma^2}\(1-\frac{a}{\sigma^2}\)\frac{1}{x^2}.
\end{equation}
whereas the Schr\"odinger potential, given by
\begin{equation}
\tilde U(z) = \(\frac{3}{4}
-\frac{4a}{\sigma^2}\(1-\frac{a}{\sigma^2}\)\)\frac{1}{z^2}
-\xi-\frac{ab}{\sigma^2}+ \frac{b^2+2\vartheta\sigma^2}{16} z^2
\end{equation}
corresponds to a discrete spectrum, following the same reasoning as
for $U(z)$. The solution $\Phi_n(z)$ is given in terms of the
Laguerre polynomials
\begin{equation}
\Phi_n(Z(x)) =\sqrt\frac{n!}{\Gamma(n+\frac{2a}{\sigma^2})}
\(\frac{\sigma^2}{2x}\)^{1/4} \(\frac{2\b
x}{\sigma^2}\)^{a/\sigma^2} e^{-\frac{\b}{\sigma^2}x}
L_n^{(\frac{2a}{\sigma^2}-1)}\(\frac{2\b}{\sigma^2}x\)
\end{equation}
with $\b=\sqrt{2\vartheta\sigma^2+b^2}$ and for the corresponding
eigenvalues $\lambda_n = -\b n - \frac{a}{\sigma^2}(\b-b)$. The
Laplace transform $L$ in (\ref{eq:Genf}) is easily integrated and
yields a convergent series in terms of the Laguerre polynomials
which re-sums to the famous formula (\ref{eq:CIRG}) (cf.
\cite{PBM1986}, p. 462(3) and p. 705(7)).
\end{proof}

\newpage
\section{Processes related to continuous orthogonal polynomials}
\label{sec:ContPoly}

\subsection{Proof of the Second Classification Theorem}
\label{sub:ProofCont}

We give a constructive proof of Theorem \ref{thm:ContPoly},
independent of Theorem \ref{thm:Diffusions}. In the following
remark, we show that Theorem \ref{thm:ContPoly} can actually be
regarded as a corollary of Theorem \ref{thm:Diffusions}.

The reducibility condition implies that the infinitesimal generator
$\L$ must be of the form
\begin{equation}
\L = \frac{\sigma^2}{2} A(x)\partII{}{x} +(a-bx)\partI{}{x}
\end{equation}
for coefficients $a\in\R$ and $b>0$. The transitional probability
density of the process, $p_{T-t}(x,y)$, satisfies the backward
Kolmogorov equation:
\begin{equation}
\partI{p}{t} +  (a-b x)\partI{p}{x} +
\frac{\sigma^2}{2} A(x)\partII{p}{x} = 0
\end{equation}
with final time condition
\begin{equation}\label{eq:ProofContFCU}
\lim_{t\to T} p_{T-t}(x,y)= \delta(x-y).
\end{equation}
By the reducibility assumption, a general solution to this equation
is given by the following eigenfunction expansion in terms of the
orthogonal polynomials $Q_n(x;a,b)$
\begin{equation}\label{eq:ProofContU}
p = \sum_{n=0}^\infty h_n(t) Q_n(x;a,b).
\end{equation}
According to (\ref{eq:EV}), the functions of time $h_n(t)$ satisfy
the ordinary differential equations
\begin{equation}
\dot h_n + \lambda_n h_n = 0
\end{equation}
which admits the general solution, for $t \le T$,
\begin{equation}
h_n(t) = z_n\ e^{\lambda_n(T-t)}.
\end{equation}
The coefficients $z_n$ are given by the final time condition
(\ref{eq:ProofContFCU}):
\begin{equation}
\sum_{n=0}^\infty z_n Q_n(x;a,b) = \delta(x-y).
\end{equation}
Hence, multiplying on both sides by $Q_m(x;a,b)$ and the invariant
measure $\rho(dx)$, before integrating over the domain of $x$, leads
to the result (\ref{eq:ContPolyU}) by the orthogonality property
(\ref{eq:Ortho}).\\

The Laplace transform $L_{T-t}(x,\vartheta)$ satisfies the backward
Kolmogorov equation with differential operator given by
(\ref{eq:G}):
\begin{equation}
\partI{L}{t} +  (a-b x)\partI{L}{x} +
\frac{\sigma^2}{2} A(x)\partII{L}{x} = \vartheta\phi(x) L
\end{equation}
with final time condition
\begin{equation}\label{eq:ProofContFC}
\lim_{t\to T} L_{T-t}(x,\vartheta)= q(x).
\end{equation}
For the sake of having a clearer constructive proof, consider the
following ansatz for the Laplace transform
\begin{equation}
L = V(x) \bar L.
\end{equation}
The previous equation then reads
\begin{equation}
V {\partial \bar L\over\partial t} + (a-b x)(V' \bar L + V \bar L')
+\frac{\sigma^2}{2}\vert A(x)\vert (V''\bar L + 2V'\bar L' + V\bar
L'') = \vartheta\phi V \bar L,
\end{equation}
where the symbol $'$ denotes differentiation in the $x$-variable.
The function $V$ is chosen to satisfy
\begin{equation}
(a-b x)V + 2\frac{\sigma^2}{2} A(x) V' = (\a-\b x)V
\end{equation}
for some parameters $\a\in\R$ and $\b>0$, which implies
\begin{equation}\label{eq:ProofContPsi}
V(x) = \exp\(\int^x\frac{(\a-\b y)-(a-b y)}{\sigma^2\vert
A(y)\vert}dy\).
\end{equation}
The function $\phi(x)$ is specified as follows:
\begin{equation}
\vartheta\phi(x) = (a-b x){V'\over V} + \frac{\sigma^2}{2}
A(x){V''\over V},
\end{equation}
which is equivalent to (\ref{eq:ContPolyphi}) with the assumption
that $C$ regroups all the constant terms. This choice yields the
following partial differential equation for the function $\bar L$:
\begin{equation}
{\partial\bar L\over\partial t} + (\a-\b x){\partial \bar
L\over\partial x} + \frac{\sigma^2}{2} A(x){\partial^2 \bar
L\over\partial x^2} = 0.
\end{equation}
Following the same reasoning as for the transitional probability
denstiy, a general solution to this equation is given by the
following eigenfunction expansion in terms of the same orthogonal
polynomials $Q_n(x;\a,\b)$, but with different coefficients,
\begin{equation}\label{eq:ProofContGb}
\bar L = \sum_{n=0}^\infty h_n(t) Q_n(x;\a,\b).
\end{equation}
According to (\ref{eq:EV}), the functions of time $h_n(t)$ satisfy
the ordinary differential equations
\begin{equation}
\dot h_n + \bar\lambda_n h_n = 0
\end{equation}
where $\bar\lambda_n$ are the eigenvalues corresponding to
$Q_n(x;\a,\b)$. The general solution, for $t \le T$, is
\begin{equation}
h_n(t) = z_n\ e^{\bar\lambda_n(T-t)}
\end{equation}
where the $z_n$ are constants. The latter equation for $h_n(t)$ with
the explicit form of $V(x)$ in (\ref{eq:ProofContPsi}) and the
expression of $\bar L$ in (\ref{eq:ProofContGb}) gives the expected
result (\ref{eq:ContPolyG}) for the Laplace transform.

The coefficients $z_n$ are given by the final time condition
(\ref{eq:ProofContFC}):
\begin{equation}
\sum_{n=0}^\infty z_n Q_n(x;\a,\b) = q(x)\ \exp\(-\int^x
\frac{(\a-\b y)-(a-b y)}{\sigma^2\vert A(y)\vert}dy\).
\end{equation}
Hence, multiplying on both sides by $Q_m(x)$ and the invariant
measure $\rho(dx)$, before integrating over the domain of $x$, leads
to the final result (\ref{eq:ContPolyzn}) by the orthogonality
property (\ref{eq:Ortho}) and concludes the proof of Theorem
\ref{thm:ContPoly}.

\begin{remark}
The reducibility condition implies that the operators $\L$ and $\G$
have the form
\begin{eqnarray}
\L &=& \frac{\sigma^2}{2}A(x)\partII{}{x}
+(a-bx)\partI{}{x}\nonumber\\
\G &=& \frac{\sigma^2}{2}A(x)\partII{}{x} +(\a-\b x)\partI{}{x}
\end{eqnarray}
for possibly different coefficients $a,\a\in\R$ and $b,\b>0$. Hence,
setting
\begin{equation}
R(x) = \frac{2A(x)}{\sigma^2}\quad{\rm and}\quad h(x)
=\exp\(\int^x\frac{a-by}{\sigma^2 A(y)}dy\)
\end{equation}
in the First Classification Theorem \ref{thm:Diffusions} yields the
Second Classification Theorem \ref{thm:ContPoly}. Moreover, consider
the gauge transformation $T_{\bar h}$ defined by
\begin{equation}
\bar h(x) =\exp\(\int^x\frac{\a-\b y}{\sigma^2 A(y)}dy\).
\end{equation}
Then we have the following relation:
\begin{equation}
T_h T_{\bar h}^{-1}\G = \frac{\sigma^2}{2}A(x)\partII{}{x} +(a-b
x)\partI{}{x} + \frac{\sigma^2}{2}\frac{Q(x,\vartheta)}{\vartheta
A(x)}.
\end{equation}
Notice further that the function $V(x)$ of the previous proof is
related to the latter two gauge transformations by $V(x)=\frac{\bar
h(x)}{h(x)}$.

\end{remark}

\subsection{The Ornstein-Uhlenbeck process}

\begin{corollary}
Assume that $\L$ and $\G$ are reducible in the sense of Definition
\ref{def:OP} to Hermite polynomials. Assume also that the function
$\phi(x)=x$. Then the transitional probability density of the
Ornstein-Uhlenbeck process is given by:
\begin{equation}\label{eq:HermiteU}
p_{T-t}(x,y) = \sqrt\frac{b}{\sigma^2\pi}e^{-\frac{b}{\sigma^2}
\(y-\frac{a}{b}\)^2}\sum_{n=0}^\infty \frac{e^{-bn(T-t)}}{n!2^n}
H_n\big(z(x)\big) H_n\big(z(y)\big)
\end{equation}
and the Laplace transform is given by the following convergent
series:
\begin{equation}\label{eq:HermiteG}
L_{T-t}(x,\vartheta) = e^{\frac{\vartheta}{b}x} \sum_{n=0}^\infty
e^{-bn(T-t)} z_n H_n\big(\bar z(x)\big),
\end{equation}
where $z(x)=\sqrt\frac{b}{\sigma^2}\big(x-\frac{a}{b}\big)$ and
$\bar z(x)=\sqrt\frac{b}{\sigma^2}\big(x-\frac{\a}{b}\big)$. The
coefficients $z_n$ are given by:
\begin{equation}\label{eq:Hermitezn}
z_n = \frac{1}{n!2^n}\sqrt\frac{b}{\sigma^2\pi}\
\int_{-\infty}^\infty q(x)\ e^{-\frac{\vartheta}{b}x} H_n\big(\bar
z(x)\big) e^{-\frac{b}{\sigma^2}\(x-\frac{\a}{b}\)^2} dx.
\end{equation}
\end{corollary}

\begin{proof}
This case fits the classification scheme in Theorem
\ref{thm:ContPoly} if one selects $A(x)=1$. This choice implies
\begin{equation}
\phi(x) = C + \frac{1}{2\vartheta\sigma^2}\big((\a-\b
x)^2-(a-bx)^2\big).
\end{equation}
In order to reduce $\phi(x)$ to an affine function, we choose
$\b=b$, $\a=a-\frac{\vartheta\sigma^2}{b}$ and
$C=\frac{a^2-\a^2}{2\vartheta\sigma^2}$. The operator $\L$, which
has the form
\begin{equation}
\L = \frac{\sigma^2}{2}\partII{}{x} + (a-bx)\partI{}{x},
\end{equation}
has the Hermite polynomials $\displaystyle H_n\big(z(x)\big)$ as
eigenfunctions with eigenvalues $\lambda_n = -bn$. The invariant
measure density and the normalization factor are respectively
\begin{equation}
w(x) = \sqrt\frac{b}{\sigma^2\pi}\
e^{-\frac{b}{\sigma^2}\(x-\frac{a}{b}\)^2}\quad{\rm and}\quad d_n^2
= n!2^n
\end{equation}
which leads, by Theorem \ref{thm:ContPoly}, to the formulation of
the kernel of the semigroup generated by $\L$ as the convergent
series in (\ref{eq:HermiteU}). The latter series re-sums using
Mehler's formula to give (\ref{eq:OUU}). By Theorem
\ref{thm:ContPoly}, the Laplace transform is given by
(\ref{eq:HermiteG}). The coefficients $z_n$ are given by
(\ref{eq:Hermitezn}) and lead to the result (\ref{eq:OUGaffine}),
once integrated and re-summed.
\end{proof}

\subsection{The CIR process}

\begin{corollary}
Assume that $\L$ and $\G$ are reducible in the sense of Definition
\ref{def:OP} to Laguerre polynomials. Assume also that the function
$\phi(x)=x$. Then the transitional probability density of the CIR
process is given by:
\begin{equation}\label{eq:LaguerreU}
p_{T-t}(x,y) = \(\frac{2b}{\sigma^2}\)^{\frac{2a}{\sigma^2}}\
y^{\frac{2a}{\sigma^2}-1}\ e^{-\frac{2b}{\sigma^2}y}
\sum_{n=0}^\infty
\frac{n!e^{-bn(T-t)}}{\Gamma(n+\frac{2a}{\sigma^2})}
L_n^{(\frac{2a}{\sigma^2}-1)}\big(\frac{2b}{\sigma^2}x\big)
L_n^{(\frac{2a}{\sigma^2}-1)}\big(\frac{2b}{\sigma^2}y\big)
\end{equation}
and the Laplace transform is given by the following convergent
series:
\begin{equation}\label{eq:LaguerreG}
L_{T-t}(x,\vartheta) = e^{\frac{b-\b}{\sigma^2}x} \sum_{n=0}^\infty
e^{-\b n(T-t)} z_n
L_n^{(\frac{2a}{\sigma^2}-1)}\big(\frac{2\b}{\sigma^2}x\big).
\end{equation}
The coefficients $z_n$ are given by:
\begin{equation}\label{eq:Laguerrezn}
z_n = \frac{n!}{\Gamma(n+\frac{2a}{\sigma^2})}
\(\frac{2\b}{\sigma^2}\)^{\frac{2a}{\sigma^2}} \int_0^\infty q(x)\
x^{\frac{2a}{\sigma^2}-1}\ e^{-\frac{b+\b}{\sigma^2}x}
L_n^{(\frac{2a}{\sigma^2}-1)}\big(\frac{2\b}{\sigma^2}x\big) dx.
\end{equation}
\end{corollary}

\begin{proof}
This case also fits the classification scheme in Theorem
\ref{thm:ContPoly} if one selects $A(x)=x$. This choice implies
\begin{equation}
\phi(x) = C + \frac{1}{2\vartheta x}\big((a-b x)-(\a-\b x)\big)+
\frac{1}{2\vartheta\sigma^2x}\big((\a-\b x)^2-(a-bx)^2\big).
\end{equation}
In order to reduce $\phi(x)$ to an affine function, we choose
$\a=a$, $\b=\sqrt{2\vartheta\sigma^2+b^2}$ and
$C=\frac{\b-b}{\vartheta}\(\frac{a}{\sigma^2}-\frac{1}{2}\)$. The
operator $\L$, which has the form
\begin{equation}
\L = \frac{\sigma^2x}{2}\ \partII{}{x} + (a-bx)\partI{}{x},
\end{equation}
has the Laguerre polynomials $\displaystyle
L_n^{(\frac{2a}{\sigma^2}-1)}\big(\frac{2b}{\sigma^2}x\big)$ as
eigenfunctions with eigenvalues $\lambda_n = -bn$ if $a>0$. The
invariant measure density and the normalization factor are
respectively
\begin{equation}
w(x) = \(\frac{2b}{\sigma^2}\)^{\frac{2a}{\sigma^2}}\
x^{\frac{2a}{\sigma^2}-1}\ e^{-\frac{2b}{\sigma^2}x}\quad{\rm
and}\quad d_n^2 = \frac{\Gamma(n+\frac{2a}{\sigma^2})}{n!}
\end{equation}
which leads, by Theorem \ref{thm:ContPoly}, to the formulation of
the kernel of the semigroup generated by $\L$ as the convergent
series in (\ref{eq:LaguerreU}) which re-sums to give
(\ref{eq:CIRU}). The Laplace transform is given by
(\ref{eq:LaguerreG}) and the coefficients $z_n$ are given by
(\ref{eq:Laguerrezn}). The latter results lead to (\ref{eq:CIRG}),
once integrated and re-summed.
\end{proof}

\subsection{The Jacobi process}
\begin{definition}\label{def:Jacobi}
The Jacobi polynomials $P_n^{(\alpha, \beta)}(x)$ are defined by the
following Gaussian hypergeometric function for $x\in[-1,1]$:
\begin{equation}
P_n^{(\alpha, \beta)}(x) = \frac{(\alpha+1)_n}{n!}{_2F_1} \Bigg(
\begin{matrix}
-n,  n+\alpha+\beta+1 \\
\alpha+1
\end{matrix}\
\Bigg\lvert\ \frac{1-x}{2} \Bigg), \ n = 0,1,2,\ldots
\end{equation}
\end{definition}

\begin{definition}
The Jacobi process is solution to the following equation:
\begin{equation}
dX_t = (a-bX_t)dt + \sigma\sqrt{X_t(1-X_t)}dW_t
\end{equation}
with initial condition $X_{t=0}=x_0\in(0,1)$.
\end{definition}

\begin{corollary}\label{cor:Jacobi}
Assume that $\L$ and $\G$ are reducible in the sense of Definition
\ref{def:OP} to Jacobi polynomials. Then:
\begin{equation}\label{eq:Jacobiphi}
\phi(x)=\frac{\sigma^2}{8\vartheta}
\(\frac{\bar\alpha^2-\alpha^2}{x}+
\frac{\bar\beta^2-\beta^2}{1-x}\),
\end{equation}
for $\alpha=\frac{2a}{\sigma^2}-1>-1,
\beta=\frac{2}{\sigma^2}(b-a)-1>-1$ and
$\bar\alpha=\frac{2\a}{\sigma^2}-1>-1,
\bar\beta=\frac{2}{\sigma^2}(\b-\a)-1>-1$. The transitional
probability density of the Jacobi process is given by:
\begin{equation}\label{eq:JacobiU}
p_{T-t}(x,y) = y^{\alpha}\ (1-y)^\beta \sum_{n=0}^\infty
\frac{e^{-\frac{\sigma^2}{2}n(n+\alpha+\beta+1)(T-t)}}{d_n^2}
P_n^{(\alpha,\beta)}(1-2x) P_n^{(\alpha,\beta)}(1-2y)
\end{equation}
with normalization constant
\begin{equation}
d_n^2 = \frac{\Gamma(n+\alpha+1)\Gamma(n+\beta+1)}
{(2n+\alpha+\beta+1)\Gamma(n+\alpha+\beta+1)n!}.
\end{equation}
The Laplace transform is given by the following convergent series:
\begin{equation}\label{eq:JacobiG}
L_{T-t}(x,\vartheta) = x^{\frac{\bar\alpha-\alpha}{2}}
(1-x)^{\frac{\bar\beta-\beta}{2}}\ \sum_{n=0}^\infty
e^{-\frac{\sigma^2}{2}n(n+\bar\alpha+\bar\beta+1)(T-t)} z_n
P_n^{(\bar\alpha,\bar\beta)}(1-2x).
\end{equation}
The coefficients $z_n$ are given by:
\begin{equation}\label{eq:Jacobizn}
z_n = \frac{1}{\bar d_n^2} \int_0^1  q(x)\
P_n^{(\bar\alpha,\bar\beta)}(1-2x)\ x^{\frac{\alpha+\bar\alpha}{2}}
(1-x)^{\frac{\beta+\bar\beta}{2}} dx,
\end{equation}
with normalization constant
\begin{equation}
\bar d_n^2 = \frac{\Gamma(n+\bar\alpha+1)\Gamma(n+\bar\beta+1)}
{(2n+\bar\alpha+\bar\beta+1)\Gamma(n+\bar\alpha+\bar\beta+1)n!}.
\end{equation}
\end{corollary}

\begin{proof}
This case also fits the classification scheme in Theorem
\ref{thm:ContPoly} if one selects $A(x)=x(1-x)$. This choice implies
\begin{equation}
\phi(x) = C + \frac{\sigma^2}{8\vartheta}
\(\frac{\bar\alpha^2-\alpha^2}{x}+ \frac{\bar\beta^2-\beta^2}{1-x}+
(\alpha+\beta)^2-(\bar\alpha+\bar\beta)^2\).
\end{equation}
Since $C$ is an arbitrary constant, we set it to
$C=(\bar\alpha+\bar\beta)^2-(\alpha+\beta)^2$. The infinitesimal
generator $\L$, which has the form
\begin{equation}
\L = \frac{\sigma^2}{2}\ x(1-x)\ \partII{}{x} + (a-bx)\partI{}{x},
\end{equation}
has the Jacobi polynomials $P_n^{(\alpha,\beta)}(1-2x)$ as
eigenfunctions with eigenvalues $\lambda_n =
-\frac{\sigma^2}{2}n(n+\alpha+\beta+1)$ if $\alpha>-1$ and
$\beta>-1$. The invariant measure density and the normalization
factor are respectively
\begin{equation}
w(x) = x^\alpha(1-x)^\beta\quad{\rm and}\quad d_n^2 =
\frac{\Gamma(n+\alpha+1)\Gamma(n+\beta+1)}
{(2n+\alpha+\beta+1)\Gamma(n+\alpha+\beta+1)n!},
\end{equation}
which concludes the proof by Theorem \ref{thm:ContPoly}.
\end{proof}

\subsection{The Dual Jacobi process}
We introduce the dual Jacobi polynomials by applying the
transformation $x\mapsto Z(x)=x(2-x)$ to the Jacobi polynomials
defined in the previous subsection.
\begin{definition}\label{def:dualJacobi}
The dual Jacobi polynomials $D_n^{(\alpha, \beta)}(x)\equiv
P_n^{(\alpha, \beta)}(1-2Z(x))$ are defined as follows:
\begin{equation}
D_n^{(\alpha, \beta)}(x) = \frac{(\alpha+1)_n}{n!}{_2F_1} \Bigg(
\begin{matrix}
-n,  n+\alpha+\beta+1 \\
\alpha+1
\end{matrix}\
\Bigg\lvert\ x(2-x) \Bigg), \ n = 0,1,2,\ldots
\end{equation}
\end{definition}
They also satisfy an orthogonality relation (\ref{eq:Ortho}) on
$(0,1)$ with normalization constants and continuous measure density:
\begin{equation}
w(x)=2\big(x(2-x)\big)^\alpha(1-x)^{2\beta+1}\quad,\quad d_n^2 =
\frac{\Gamma(n+\alpha+1)\Gamma(n+\beta+1)}
{(2n+\alpha+\beta+1)\Gamma(n+\alpha+\beta+1)n!}.
\end{equation}
They are solutions to the eigenvalue problem (\ref{eq:EV}) with
generator:
\begin{equation}\label{eq:dualJacobiL}
\L = \frac{\sigma^2}{8}x(2-x)\partII{}{x}
+\frac{2a-(2b-\frac{\sigma^2}{2})x(2-x)}{4(1-x)}\partI{}{x}
\end{equation}
and eigenvalues $\lambda_n=-\frac{\sigma^2}{2}n(n+\alpha+\beta+1)$,
still conditioned to $\alpha=\frac{2a}{\sigma^2}-1>-1$ and
$\beta=\frac{2}{\sigma^2}(b-a)-1>-1$. Hence, we have the definition
of the dual Jacobi process as follows:
\begin{definition}
The dual Jacobi process is solution to the following equation:
\begin{equation}
dX_t = \frac{2a-(2b-\frac{\sigma^2}{2})X_t(2-X_t)}{4(1-X_t)}dt +
\frac{\sigma}{2} \sqrt{X_t(2-X_t)}dW_t
\end{equation}
with initial condition $X_{t=0}=x_0\in(0,1)$.
\end{definition}
We obtain yet another corollary to Theorem \ref{thm:ContPoly}:

\begin{corollary}\label{cor:dualJacobi}
Assume that $\L$ and $\G$ are reducible in the sense of Definition
\ref{def:OP} to the dual Jacobi polynomials. Then:
\begin{equation}\label{eq:dualJacobiphi}
\phi(x)=\frac{\sigma^2}{8\vartheta}
\(\frac{\bar\alpha^2-\alpha^2}{x(2-x)}+
\frac{\bar\beta^2-\beta^2}{(1-x)^2}\),
\end{equation}
for $\bar\alpha=\frac{2\a}{\sigma^2}-1>-1,
\bar\beta=\frac{2}{\sigma^2}(\b-\a)-1>-1$. The transitional
probability density of the dual Jacobi process is given by:
\begin{equation}\label{eq:dualJacobiU}
p_{T-t}(x,y) = 2\big(y(2-y)\big)^{\alpha}\ (1-y)^{2\beta+1}
\sum_{n=0}^\infty \frac{e^{\lambda_n(T-t)}}{d_n^2}
D_n^{(\alpha,\beta)}(x) D_n^{(\alpha,\beta)}(y).
\end{equation}
The Laplace transform can be expressed by the following convergent
series:
\begin{equation}\label{eq:dualJacobiG}
L_{T-t}(x,\vartheta) = 2
\big(x(2-x)\big)^{\frac{\bar\alpha-\alpha}{2}}
(1-x)^{\bar\beta-\beta}\ \sum_{n=0}^\infty e^{\lambda_n(T-t)} z_n
D_n^{(\bar\alpha,\bar\beta)}(x)
\end{equation}
where the coefficients $z_n$ are given by:
\begin{equation}\label{eq:dualJacobizn}
z_n = \frac{1}{\bar d_n^2}\int_0^1 q(x)\
D_n^{(\bar\alpha,\bar\beta)}(x)\
\big(x(2-x)\big)^{\frac{\alpha+\bar\alpha}{2}}
(1-x)^{\beta+\bar\beta+1}dx.
\end{equation}
with normalization constant
\begin{equation}
\bar d_n^2 = \frac{\Gamma(n+\bar\alpha+1)\Gamma(n+\bar\beta+1)}
{(2n+\bar\alpha+\bar\beta+1)\Gamma(n+\bar\alpha+\bar\beta+1)n!}.
\end{equation}
\end{corollary}
\begin{proof}
The proof follows from Corollary \ref{cor:Jacobi} and the
transformation $x\mapsto Z(x)=x(2-x)$.
\end{proof}

\newpage
\section{Processes related to discrete orthogonal polynomials}
\label{sec:DiscPoly}

The proof of Theorem \ref{thm:DiscPoly} follows a very similar
reasoning as the proof of Theorem \ref{thm:ContPoly}, its continuous
version.

\subsection{Proof of the Third Classification Theorem}

The reducibility condition implies that the infinitesimal generator
$\L$ must be of the form
\begin{equation}
\L = -D(x)\Delta^1 + \big(D(x)-B(x)\big)\nabla_+^1.
\end{equation}
The transitional probability density $p_{T-t}(x,y)$ satisfies the
backward Kolmogorov equation with generator given by (\ref{eq:L}):
\begin{equation}
\partI{p}{t} -D(x)\Delta^1p + \big(D(x)-B(x)\big)\nabla_+^1p = 0
\end{equation}
with final time condition
\begin{equation}\label{eq:ProofDiscFCU}
\lim_{t\to T}p_{T-t}(x,y)= \delta(x-y).
\end{equation}
By the reducibility assumption, a general solution to this equation
is given by the following eigenfunction expansion in terms of the
discrete orthogonal polynomials $Q_n(x)$
\begin{equation}\label{eq:ProofDiscU}
p = \sum_{n=0}^\infty h_n(t) Q_n(x).
\end{equation}
According to (\ref{eq:EV}), the functions of time $h_n(t)$ satisfy
the ordinary differential equations
\begin{equation}
\dot h_n + \lambda_n h_n = 0
\end{equation}
which admits the general solution, for $t \le T$,
\begin{equation}
h_n(t) = z_n\ e^{\lambda_n(T-t)}.
\end{equation}
The coefficients $z_n$ are given by the final time condition
(\ref{eq:ProofDiscFCU}):
\begin{equation}
\sum_{n=0}^\infty z_n Q_n(x) = \delta(x-y).
\end{equation}
Hence, multiplying on both sides by $Q_m(x)$ and the weight $w(x)$,
before summing over the lattice $\Lambda_N$, leads to the result
(\ref{eq:DiscPolyU}) by the orthogonality property
(\ref{eq:Ortho}).\\

The Laplace transform $L_{T-t}(x,1)$ satisfies this time a finite
difference version of the Backward Kolmogorov equation
\begin{equation}
\partI{L}{t} -D(x)\Delta^1L + \big(D(x)-B(x)\big)\nabla_+^1L =
\vartheta\phi L
\end{equation}
with the same final time condition
\begin{equation}\label{eq:ProofDiscFC}
\lim_{t\to T}L_{T-t}(x,\vartheta)= q(x).
\end{equation}
Consider the ansatz for the Laplace transform
\begin{equation}
L = V(x) \bar L.
\end{equation}
The latter finite difference equation reads
\begin{equation}
\partI{\bar L}{t}
- \underbrace{B(x)\frac{V(x+1)}{V(x)}}_{\bar B(x)}\bar L(x+1) +
\underbrace{\big(B(x)+D(x)- \vartheta\phi\big)}_{\bar B(x)+\bar
D(x)}\bar L(x) - \underbrace{D(x)\frac{V(x-1)}{V(x)}}_{\bar
D(x)}\bar L(x-1) = 0.
\end{equation}
$\bar B(x)$ and $\bar D(x)$ are defined such that they satisfy the
relations
\begin{equation}
\frac{V(x)}{V(x-1)} = \frac{\bar B(x-1)}{B(x-1)} = \frac{D(x)}{\bar
D(x)}
\end{equation}
which implies condition (\ref{eq:ConditionBD}) and solves
iteratively to give
\begin{equation}\label{eq:ProofDiscPsi}
V(x) = \prod_{k=1}^x\frac{D(k)}{\bar D(k)}.
\end{equation}
The function $\phi(x)$ is specified as follows:
\begin{equation}
\phi(x) = \frac{1}{\vartheta}\big(B(x)+ D(x)-\bar B(x) - \bar
D(x)\big).
\end{equation}
But $\phi(x), B(x), D(x)$ are all by definition independent of the
parameter $\vartheta$, so we are bound to set $\vartheta=1$. This
choice yields the following finite difference equation for the
function $\bar L$:
\begin{equation}\label{eq:ProofLGFD}
\partI{\bar L}{t} - \bar D(x)\Delta^1\bar L +
\big(\bar D(x)-\bar B(x)\big)\nabla_+^1\bar L = 0.
\end{equation}
A general solution to this equation is given by the following
eigenfunction expansion in terms of discrete orthogonal polynomials:
\begin{equation}\label{eq:ProofDiscGb}
\bar L = \sum_{n=0}^N h_n(t) \bar Q_n(x).
\end{equation}
The functions of time $h_n(t)$ satisfy the ordinary differential
equations
\begin{equation}
\dot h_n +\bar \lambda_n\ h_n =0
\end{equation}
which admits the general solution
\begin{equation}
h_n(t) = z_n\ e^{\bar \lambda_n(T-t)}
\end{equation}
where the $z_n$ are constants. The latter equation for $h_n(t)$ with
the explicit form of $V(x)$ in (\ref{eq:ProofDiscPsi}) and the
expression of $\bar L$ in (\ref{eq:ProofDiscGb}) yields the
expression (\ref{eq:DiscPolyG}) for the Laplace transform.

The coefficients $z_n$ are given by the final time condition
(\ref{eq:ProofDiscFC}):
\begin{equation}
q(x)\ \prod_{k=1}^x\frac{\bar D(x)}{D(x)} = \sum_{n=0}^N z_n \bar
Q_n(x).
\end{equation}
Finally, multiplying on both sides by $\bar Q_m(x)$ and the weight
$\bar w(x)$, before summing over $\Lambda_N$, gives the final result
(\ref{eq:DiscPolyzn}) by orthogonality of the polynomials and
concludes the proof of Theorem \ref{thm:DiscPoly}.

\subsection{The Meixner process}
The Meixner polynomials provide a discrete lattice approximation to
the Laguerre polynomials.
\begin{definition}
The Meixner polynomials are defined as follows in case $x$ is
integer:
\begin{equation}
M_n(x; \beta, c) =  {_2F_1} \Bigg(
\begin{matrix}
-n, -x \\
\beta
\end{matrix}\
\Bigg\vert\ 1-\frac{1}{c} \Bigg), \ \ n = 0, 1, 2, \ldots
\end{equation}
\end{definition}
The Meixner polynomials satisfy an orthogonality relation with
respect to the discrete measure supported on $\Z_+$. Namely,
\begin{equation}
\sum_{x=0}^\infty\ M_m(x;\beta,c)\ M_n(x; \beta, c)\  w(x)  =
{c^{-n} n! \over (\beta)_n (1-c)^\beta } \delta_{nm}.
\end{equation}
where the weight is
\begin{equation}
w(x) = {(\beta)_x\over x! } c^x.
\end{equation}
The Meixner polynomials are solutions to the eigenvalue problem
(\ref{eq:EV}) with generator
\begin{equation}
\L =\frac{\sigma^2}{2}\ x\ \Delta^1 + (a-bx)\nabla^1_+ ,
\end{equation}
for $x\in\Z_+$, with $a,b>0$ and eigenvalues $\lambda_n = -bn$. The
latter can be recast in the form (\ref{eq:DiscL}) using the
functions
\begin{eqnarray}
B(x)&=& -\frac{\sigma^2}{2}c(x+\beta)\nonumber\\
D(x)&=& -\frac{\sigma^2}{2}x,
\end{eqnarray}
for $a=\frac{\sigma^2}{2} \beta$ and $b=\frac{\sigma^2}{2} (1-c)$.
The parameters are conditioned to $\beta>0$ and $0<c<1$ which
insures the Markov property. The Meixner process is the discrete
Markov process generated by $\L$. It is a discrete version of the
CIR process.\\

The following statement is a corollary to Theorem
\ref{thm:DiscPoly}.
\begin{corollary}\label{cor:Meixner}
Assume that $\L$ and $\G$ are reducible in the sense of Definition
\ref{def:OP} to the Meixner polynomials. Assume also that the
function $\phi(x)$ is given by $\phi(x) = \varrho x+\zeta$,
\begin{equation}
\varrho = \frac{\sigma^2}{2}\big(c(e^\varphi-1) +
(e^{-\varphi}-1)\big)\quad,\quad  \zeta = \frac{\sigma^2}{2}\beta
c(e^\varphi-1)
\end{equation}
with the real parameter $\varphi<-\frac{1}{2}\ln c$. Then the
transitional probability density for the Meixner process is as
follows:
\begin{eqnarray}\label{eq:MeixnerU}
p_{T-t}(x,y) &=& (1-c)^\beta\ \frac{(1-e^{(T-t)(c-1)})^{x+y}}
{(1-ce^{(T-t)(c-1)})^{x+y+\beta}}\
\frac{(\beta)_{y} c^{y}}{y!}\nonumber\\
&&\cdot\ {_2F_1}\Bigg(\begin{matrix}-x,-y\\
\beta\end{matrix}\ \Bigg\vert\ \frac{e^{(T-t)(c-1)}(1-c)^2}{c
(1-e^{(T-t)(c-1)})^2} \Bigg).
\end{eqnarray}
For $q(x) = \exp\big(\omega \phi(x)\big)$, the Laplace transform is
affine:
\begin{equation}\label{eq:MeixnerG}
L_{T-t}(x,1,\omega) = e^{m(T-t;\omega)x + n(T-t;\omega)}.
\end{equation}
The functions of time $m(\tau;\omega)$ and $n(\tau;\omega)$ are as
follows:
\begin{eqnarray}
m(\tau;\omega) &=& \log\(e^\varphi\
\frac{1-ce^{\varrho\omega+\varphi} -
e^{(\c-1)e^{-\varphi}\tau}(1-e^{\varrho\omega-\varphi})}
{1-ce^{\varrho\omega+\varphi} - \c
e^{(\c-1)e^{-\varphi}\tau}(1-e^{\varrho\omega-\varphi})}\)\nonumber\\
n(\tau;\omega) &=& -\beta\log\(e^{\frac{\sigma^2}{2}\omega
c(1-e^\varphi)}\ \frac{1-ce^{\varrho\omega+\varphi} - \c
e^{(\c-1)e^{-\varphi}\tau}(1-e^{\varrho\omega-\varphi})}{1-\c}\).
\nonumber
\end{eqnarray}
where $\c = ce^{2\varphi}$.
\end{corollary}

\begin{proof}
The transitional probability density follows from equation
(\ref{eq:DiscPolyU}) in the discrete classification theorem. The
definition of $\phi(x)$ suggests that we set $\bar B(x) =
B(x)e^\varphi$ and $\bar D(x) = D(x) e^{-\varphi}$ in order to
satisfy condition (\ref{eq:ConditionBD}) in Theorem
\ref{thm:DiscPoly}. The generator defined by $\bar B(x)$ and $\bar
D(x)$ has eigenfunctions $M_n(x;\beta,\c)$ with eigenvalues
$\bar\lambda_n = \frac{\sigma^2}{2}e^{-\varphi}n(\c-1)$. Now from
(\ref{eq:DiscPolyzn}), we have
\begin{eqnarray}
z_n &=&\frac{(\beta)_n \c^n (1-\c)^\beta}{n!} \sum_{x=0}^\infty
e^{\omega\phi(x)-\varphi x}\ M_n(x;\beta,\c)\
\frac{(\beta)_x\c^x}{x!}\nonumber\\
&=&e^{\omega\zeta}\ \(\frac{1-\c}{1-\c
e^{\varrho\omega-\varphi}}\)^\beta\ \frac{(\beta)_n \c^n}{n!} \
\(\frac{1-e^{\varrho\omega-\varphi}} {1-\c
e^{\varrho\omega-\varphi}}\)^n.\nonumber
\end{eqnarray}
The Laplace transform, given by (\ref{eq:DiscPolyG}), is as follows:
\begin{eqnarray}
L_{T-t}(x,1,\omega) &=& \(\frac{1-\c}{1-\c
e^{\varrho\omega-\varphi}}\)^\beta\ e^{\omega\zeta}\ e^{\varphi
x}\nonumber\\
&&\cdot\ \sum_{n=0}^\infty \frac{(\beta)_n}{n!}\ \(\c
e^{\frac{\sigma^2}{2}(\c-1)e^{-\varphi}(T-t)}
\frac{1-e^{\varrho\omega-\varphi}}{1-\c
e^{\varrho\omega-\varphi}}\)^n \
M_n(x;\beta,\c)\nonumber\\
&=& \(e^{\frac{\sigma^2}{2}\omega c(1-e^\varphi)}\
\frac{1-ce^{\varrho\omega+\varphi} - \c
(1-e^{\varrho\omega-\varphi}) e^{\frac{\sigma^2}{2}
(\c-1)e^{-\varphi}(T-t)}}{1-\c}\)^{-\beta} \nonumber\\ && \cdot\
\(e^\varphi\ \frac{1-ce^{\varrho\omega+\varphi} -
(1-e^{\varrho\omega-\varphi})
e^{\frac{\sigma^2}{2}(\c-1)e^{-\varphi}(T-t)}}
{1-ce^{\varrho\omega+\varphi} - \c (1-e^{\varrho\omega-\varphi})
e^{\frac{\sigma^2}{2}(\c-1)e^{-\varphi}(T-t)} }\)^{x}.\nonumber
\end{eqnarray}
Note that the re-summation formula used to find the last two results
is the generating function for the Meixner polynomials which can be
found in \cite{KS1998}. Also notice that
$M_n(x;\beta,\c)=M_x(n;\beta,\c)$ by definition.
\end{proof}

\subsection{The Racah Process}
\begin{definition}\label{def:Racah}
The Racah polynomials $R_n(\lambda(x)):=R_n(\lambda(x);\alpha, \beta,
\gamma, \delta)$ are defined as follows:
\begin{equation}
R_n(\lambda(x);\alpha,\beta,\gamma,\delta) = {_4F_3} \left(
\begin{matrix}
-n,  n+\alpha+\beta+1, -x, x+\gamma+\delta+1 \\
\alpha+1, \beta+\delta+1, \gamma+1
\end{matrix}
\bigg\lvert 1 \right), \ n = 0,1,2,\ldots,N
\end{equation}
where $\lambda(x) = x(x+\gamma+\delta+1)$ and either $\alpha = -N-1$
or $\beta+\delta = -N - 1$ or $\gamma = -N-1$.
\end{definition}
The Racah polynomials satisfy an orthogonality relation with respect to
the discrete measure supported on the set $\Lambda_N$. Namely,
\begin{equation}
\sum_{x\in\Lambda_N} R_m(\lambda(x)) R_n(\lambda(x)) w(x) = d_n^2
\delta_{nm}
\end{equation}
where the weight is
\begin{equation}
w(x):=w(x;\alpha,\beta,\gamma,\delta) = \frac{(\alpha+1)_x
(\beta+\delta+1)_x (\gamma+1)_x (\gamma+\delta+1)_x
((\gamma+\delta+3)/2)_x}{(-\alpha+\gamma+\delta+1)_x(-\beta+\gamma+1)_x
((\gamma+\delta+1)/2)_x (\delta+1)_x x!}
\end{equation}
and the normalization factor is
\begin{equation}
d_n^2 = M\ \frac{(n+\alpha+\beta+1)_n(\alpha+\beta-\gamma+1)_n
(\alpha-\delta+1)_n (\beta+1)_n n!} {(\alpha+\beta+2)_{2n}
(\alpha+1)_n (\beta+\delta+1)_n (\gamma+1)_n}
\end{equation}
with
$$
M = \left\{
\begin{array}{ll}
\displaystyle\frac{(-\beta)_N(\gamma+\delta+2)_N}
{(-\beta+\gamma+1)_N(\delta+1)_N} &
{\rm if}\ \alpha=-N-1\\
\displaystyle\frac{(-\alpha+\delta)_N(\gamma+\delta+2)_N}
{(-\alpha+\gamma+\delta+1)_N(\delta+1)_N} &
{\rm if}\ \beta+\delta=-N-1\\
\displaystyle\frac{(\alpha+\beta+2)_N(-\delta)_N}
{(\alpha-\delta+1)_N(\beta+1)_N} &
{\rm if}\ \gamma=-N-1.\\
\end{array}
\right.
$$
The Racah polynomials are solutions to the eigenvalue problem
(\ref{eq:EV}) with generator (\ref{eq:DiscL}) given by the functions
\begin{eqnarray}
B(x)&=&
\frac{\sigma^2}{2}\frac{(x+\alpha+1)(x+\beta+\delta+1)(x+\gamma+1)
(x+\gamma+\delta+1)} {(2x+\gamma+\delta+1)(2x+\gamma+\delta+2)}, \nonumber\\
D(x)&=&
\frac{\sigma^2}{2}\frac{x(x-\alpha+\gamma+\delta)(x-\beta+\gamma)(x+\delta)}
{(2x+\gamma+\delta)(2x+\gamma+\delta+1)},
\end{eqnarray}
and eigenvalues $\lambda_n =-\frac{\sigma^2}{2}n(n+\alpha+\beta+1)$.
The Markov property is insured if $B(x)\le0$ and $D(x)\le0$,
$\forall x\in\Lambda_N$. The process generated by the latter
generator is called the {\it Racah process}. Also define the
corresponding functions
\begin{eqnarray}\label{eq:RacahBD}
\bar B(x)&=&
\frac{\sigma^2}{2}\frac{(x+\bar\alpha+1)(x+\bar\beta+\bar\delta+1)
(x+\bar\gamma+1)(x+\bar\gamma+\bar\delta+1)}
{(2x+\bar\gamma+\bar\delta+1)(2x+\bar\gamma+\bar\delta+2)}, \nonumber\\
\bar D(x)&=&
\frac{\sigma^2}{2}\frac{x(x-\bar\alpha+\bar\gamma+\bar\delta)
(x-\bar\beta+\bar\gamma)(x+\bar\delta)}
{(2x+\bar\gamma+\bar\delta)(2x+\bar\gamma+\bar\delta+1)}
\end{eqnarray}
for $\bar\alpha,\bar\beta,\bar\gamma,\bar\delta\in\R$.

\begin{definition}
The set of parameters
$\{\bar\alpha,\bar\beta,\bar\gamma,\bar\delta\}$ will be called {\it
acceptable} with respect to $\{\alpha,\beta,\gamma,\delta\}$ if it
satisfies condition (\ref{eq:ConditionBD}) in Theorem
\ref{thm:DiscPoly}, i.e. if $\bar B(x-1)\bar D(x) = B(x-1)D(x)$, and
if $\bar B(x)\le0$ and $\bar D(x)\le0$, $\forall x\in\Lambda_N$.
\end{definition}
The following statement is a corollary to Theorem
\ref{thm:DiscPoly}.
\begin{corollary}\label{cor:Racah}
Assume that $\L$ and $\G$ are reducible in the sense of Definition
\ref{def:OP} to Racah polynomials. Assume also that for an
acceptable set of parameters
$\{\bar\alpha,\bar\beta,\bar\gamma,\bar\delta\}$:
\begin{equation}\label{eq:Racahphi}
\phi(x) = B(x)+D(x)-\bar B(x)-\bar D(x).
\end{equation}
Then the transitional probability density for the Racah process is
given by
\begin{equation}\label{eq:RacahU}
p_{T-t}(x,y) =\sum_{n=0}^N
\frac{e^{-\frac{\sigma^2}{2}n(n+\alpha+\beta+1)(T-t)}}{d_n^2}
R_n(\lambda(x)) R_n(\lambda(y)) w(y).
\end{equation}
while the Laplace transform is given by the following convergent
series:
\begin{equation}\label{eq:RacahG}
L_{T-t}(x,1) =\prod_{k=1}^{x}\frac{D(k)}{\bar D(k)}\ \sum_{n=0}^N
e^{-\frac{\sigma^2}{2}n(n+\bar\alpha+\bar\beta+1)(T-t)} z_n
R_n(\bar\lambda(x);\bar\alpha,\bar\beta,\bar\gamma,\bar\delta).
\end{equation}
The coefficients $z_n$ are as follows:
\begin{equation}\label{eq:Racahzn}
z_n = \frac{1}{\bar d_n^2}\sum_{x\in\Lambda_N}
\prod_{k=1}^x\frac{\bar D(k)}{D(k)}\ q(x)
R_n(\bar\lambda(x);\bar\alpha,\bar\beta,\bar\gamma,\bar\delta) \bar
w(x)
\end{equation}
where $\bar d_n = d_n(\bar\alpha,\bar\beta,\bar\gamma,\bar\delta)$,
$\bar w(x) = w(x;\bar\alpha,\bar\beta,\bar\gamma,\bar\delta)$ and
$\bar\lambda(x)=\lambda(x;\bar\gamma,\bar\delta)$.
\end{corollary}

\begin{proof}
The restrictions imposed on the set of parameters
$\{\bar\alpha,\bar\beta,\bar\gamma,\bar\delta\}$ ensures that
condition (\ref{eq:ConditionBD}) in Theorem \ref{thm:DiscPoly} is
satisfied. This is a necessary condition. The rest of the corollary
is a direct application of Theorem \ref{thm:DiscPoly}.
\end{proof}

\subsection{The Dual Hahn Process}

\begin{definition}\label{def:dualHahn}
The dual Hahn polynomials $R_n(\lambda(x)):=R_n(\lambda(x); \gamma,
\delta,N)$ are defined as follows:
\begin{equation}
R_n(\lambda(x);\gamma,\delta,N) = {_3F_2} \(
\begin{matrix}
-n, -x, x+\gamma+\delta+1 \\
\gamma+1, -N
\end{matrix}
\bigg\lvert 1 \), \ n = 0,1,2,\ldots,N
\end{equation}
where $\lambda(x) = x(x+\gamma+\delta+1)$.
\end{definition}

For $\gamma>-1$ and $\delta>-1$ or for $\gamma<-N$ and $\delta<-N$, the
dual Hahn polynomials satisfy an orthogonality relation with respect to
the discrete measure supported on the set $\Lambda_N$. Namely,
\begin{equation}
\sum_{x\in\Lambda_N} R_m(\lambda(x)) R_n(\lambda(x)) w(x) = d_n^2
\delta_{nm}
\end{equation}
where the weight is
\begin{equation}
w(x):=w(x;\gamma,\delta,N) = \frac{(2x+\gamma+\delta+1) (\gamma+1)_x
(-N)_x N!} {(-1)^x(x+\gamma+\delta+1)_{N+1}(\delta+1)_x x!}
\end{equation}
and the normalization factor is
\begin{equation}
d_n^2 = \frac{1}{\( \begin{matrix}\gamma+n\\n\end{matrix}\)
\(\begin{matrix}\delta+N-n\\N-n\end{matrix}\)}.
\end{equation}
The dual Hahn polynomials are solutions to the eigenvalue problem
(\ref{eq:EV}) with generator (\ref{eq:DiscL}) given by the functions
\begin{eqnarray}
B(x)&=&
-\frac{\sigma^2}{2}\frac{(x+\gamma+1)(x+\gamma+\delta+1)(N-x)}
{(2x+\gamma+\delta+1)(2x+\gamma+\delta+2)}, \nonumber\\
D(x)&=& -\frac{\sigma^2}{2}\frac{x(x+\gamma+\delta+N+1)(x+\delta)}
{(2x+\gamma+\delta)(2x+\gamma+\delta+1)},
\end{eqnarray}
and eigenvalues $\lambda_n=-\frac{\sigma^2}{2}n$. The discrete
Markov process generated by the latter generator is called the {\it
dual Hahn process}.

The following statement is another corollary to Theorem
\ref{thm:DiscPoly}.
\begin{corollary}\label{cor:dualHahn}
Assume that $\L$ and $\G$ are reducible in the sense of Definition
\ref{def:OP} to dual Hahn polynomials. Let $\delta>\gamma>-1$ or
$\delta<\gamma<-N$ such that:
\begin{equation}
\phi(x) =
\frac{\sigma^2}{2}(\delta-\gamma)\[\frac{x(x+\gamma+\delta+N+1)}
{(2x+\gamma+\delta)(2x+\gamma+\delta+1)} -
\frac{(x+\gamma+\delta+1)(N-x)}
{(2x+\gamma+\delta+1)(2x+\gamma+\delta+2)}\].
\end{equation}
Then the transitional probability density for the dual Hahn process
is given by
\begin{equation}\label{eq:dualHahnU}
p_{T-t}(x,y) =\sum_{n=0}^N
\frac{e^{-\frac{\sigma^2}{2}n(T-t)}}{d_n^2} R_n(\lambda(x))
R_n(\lambda(y)) w(y).
\end{equation}
while the Laplace transform can be expressed as the following
convergent series:
\begin{equation}\label{eq:dualHahnG}
L_{T-t}(x,1) =\frac{(\delta+1)_{x}}{(\gamma+1)_{x}}\ \sum_{n=0}^N
e^{-\frac{\sigma^2}{2}n(T-t)} z_n R_n(\lambda(x);\delta,\gamma,N).
\end{equation}
The coefficients $z_n$ are as follows:
\begin{equation}\label{eq:dualHahnzn}
z_n = \frac{1}{\bar d_n^2}\sum_{x\in\Lambda_N}
\frac{(\gamma+1)_x}{(\delta+1)_x}\ q(x)
R_n(\lambda(x);\delta,\gamma,N) w(x;\delta,\gamma,N)
\end{equation}
where $\bar d_n = d_n(\delta,\gamma,N)$.
\end{corollary}

\begin{proof}
For $\bar\gamma=\delta$ and $\bar\delta=\gamma$, define the
functions
\begin{eqnarray}
\bar B(x)&=&
-\frac{\sigma^2}{2}\frac{(x+\bar\gamma+1)(x+\bar\gamma+\bar\delta+1)
(N-x)}{(2x+\bar\gamma+\bar\delta+1)
(2x+\bar\gamma+\bar\delta+2)}, \nonumber\\
\bar D(x)&=&
-\frac{\sigma^2}{2}\frac{x(x+\bar\gamma+\bar\delta+N+1)(x+\bar\delta)}
{(2x+\bar\gamma+\bar\delta)(2x+\bar\gamma+\bar\delta+1)}.
\end{eqnarray}
The corollary then follows from Theorem \ref{thm:DiscPoly}.
\end{proof}

\newpage
\section{Limit relations}\label{sec:Limits}
Limit relations between orthogonal polynomials are well-known, see
\cite{KS1998} and \cite{NSU1991}. In this section, we show that the
transitional probability densities and Laplace transforms obtained
from the various corollaries of Theorems \ref{thm:ContPoly} and
\ref{thm:DiscPoly} are similarly connected to each other. We start
by rigorously stating the relation between the Jacobi process and
its dual.

\begin{proposition}
The dual Jacobi process is obtained from the Jacobi process with the
change of variable transformation
\begin{equation}
Z(x) = x(2-x)
\end{equation}
applied to the underlying process $X_t$. The same holds for its
transitional probability density and Laplace transform.
\end{proposition}

For the discrete $X_t$ process, the limit relation between the Racah
process and the dual Hahn process is not as obvious.

\begin{proposition}
The Racah process converges to the dual Hahn process in three
different ways, corresponding to the three families of Racah
polynomials:
\begin{enumerate}
\item $\alpha=-N-1$ with the acceptable set of parameters
$\{-N-1,\beta+\delta-\gamma,\delta,\gamma\}$, conditioned to either
$$
\left\{\begin{array}{l}
\beta\ge\gamma+N\\
\delta>\gamma>-1
\end{array}\right.
\quad or\quad \left\{\begin{array}{l}
\beta\ge-\delta-1\\
\delta<\gamma<-N
\end{array}\right.
$$
in the limit as $\beta\to\infty$.

\item $\beta=-\delta-N-1$ with the acceptable set $\{\alpha, -\gamma-N-1,
\delta,\gamma\}$, conditioned to either
$$
\left\{\begin{array}{l}
\alpha\ge\gamma+\delta+N\\
\delta>\gamma>-1
\end{array}\right.
\quad or\quad \left\{\begin{array}{l}
\alpha\ge-1\\
\delta<\gamma<-N
\end{array}\right.
$$
in the limit as $\alpha\to\infty$.

\item $\gamma=-N-1$ with the acceptable set $\{-\alpha+\delta-N-1, \beta,
-N-1,\delta\}$, conditioned to either
$$
\left\{\begin{array}{l}
\alpha>-1\\
\beta\ge-1\\
\delta>2\alpha+N+1
\end{array}\right.
\quad or\quad \left\{\begin{array}{l}
\alpha<-N\\
\beta\ge-\delta-1\\
\delta<2\alpha+N+1\\
\end{array}\right.,
$$
with first the mapping $\delta\mapsto\delta+\alpha+N+1$ and then the
limit $\beta\to\infty$. The dual Hahn parameters are in this third
case $(\alpha,\delta)$.
\end{enumerate}
The result extends to the transitional probability densities and
Laplace transforms.
\end{proposition}

The next proposition states that the Meixner process converges to
the CIR process in the affine case.
\begin{proposition}\label{prop:CIRMeixner}
Under the transformations $\varrho\mapsto\varrho(1-c)$ and
$x\mapsto\displaystyle\frac{x}{1-c}$ in the Meixner process, the
limit $c\to1$ yields the CIR process with parameter
$\alpha=\beta-1$. The Laplace transform is affine in this case.
\end{proposition}
\begin{proof}
The proof follows from the limit relation
$$
\lim_{c\to1} M_n\(\frac{x}{1-c};\alpha+1,c\)=
\frac{L^{(\alpha)}_n(x)}{L^{(\alpha)}_n(0)}.
$$
Both the transitional probability density and the Laplace transform
in the Meixner case then converge to the affine CIR case, since
$$
\lim_{c\to1} \phi\(\frac{x}{1-c}\) = \varrho\ x + \zeta.
$$
\end{proof}

The last proposition demonstrates the connection between a discrete
and a continuous underlying process $X_t$ by showing that the dual
Jacobi process is actually a limiting case of the Racah process.

\begin{proposition}\label{prop:RacahdualJacobi}
Consider the acceptable set of parameters
$\{\alpha,\bar\beta,\gamma,\delta\}$ where
\begin{eqnarray}
\bar\beta&=&-\beta-\delta-N-1\nonumber\\
\gamma&=&-N-1.
\end{eqnarray}
Assume furthermore the inequalities:
\begin{eqnarray}
\alpha&>&-1\nonumber\\
\bar\beta&>&\beta\ >\ -1\nonumber\\
\gamma&>&\delta.
\end{eqnarray}
Then, applying the transformation $x\mapsto xN$ to the Racah process
such that $x\in[0,1]$, yields the dual Jacobi process in the limit
$N\to\infty$ for the special case of $\bar\alpha=\alpha$ and
$\bar\beta>\vert\beta\vert$. This result applies to both the
transitional probability density and the Laplace transform.
\end{proposition}

\begin{proof}
The inequalities in the assumption insure that the process $X_t$
satisfies the Markov property, as both $B(x)$ and $D(x)$ are
negative for $x\in\Lambda_N$. Since on top of $\bar\alpha=\alpha,
\bar\gamma=\gamma, \bar\delta=\delta$, we have $-\beta+\gamma =
\bar\beta+\bar\delta$ and $\beta+\delta=-\bar\beta+\bar\gamma$, it
is immediate that $\{\alpha,\bar\beta,\gamma,\delta\}$ is
acceptable. The function $\phi(x)$ reduces to
\begin{equation}
\phi(x) = \frac{\sigma^2}{2}(\bar\beta-\beta)\
\[\frac{x(x-\alpha+\gamma+\delta)(x+\delta)}
{(2x+\gamma+\delta)(2x+\gamma+\delta+1)}-
\frac{(x+\alpha+1)(x+\gamma+1)(x+\gamma+\delta+1)}
{(2x+\gamma+\delta+1)(2x+\gamma+\delta+2)}\]
\end{equation}
or equivalently,
\begin{eqnarray}
\phi(x) &=& \frac{\sigma^2}{2}(\bar\beta-\beta)\
\[\frac{x(x-2N-\alpha-\bar\beta-\beta-2)(x-N-\bar\beta-\beta-1)}
{(2x-2N-\bar\beta-\beta-2)(2x-2N-\bar\beta-\beta-1)}\right.\nonumber\\
&&- \left.\frac{(x+\alpha+1)(x-N)(x-2N-\bar\beta-\beta-1)}
{(2x-2N-\bar\beta-\beta-1)(2x-2N-\bar\beta-\beta)}\].
\end{eqnarray}
$\phi(x)$ is bounded from below, from the inequality
$\bar\beta>\beta>-1$ and is even monotonously increasing. Thus the
assumptions of Corollary \ref{cor:Racah} are satisfied. So we have
defined a proper Racah process with a Laplace transform that can be
expressed as a convergent series over the Racah polynomials. With
the transformation $x\mapsto Nx$, the function $\phi(Nx)$ converges,
in the limit $N\to\infty$, to its continuous counterpart of
Corollary \ref{cor:dualJacobi} for $\vartheta=1$.

Moreover, the finite difference infinitesimal generator $\L$ in
(\ref{eq:DiscL}) converges to the dual Jacobi diffusion generator in
(\ref{eq:dualJacobiL}), under the transformation $x\mapsto Nx$:
\begin{eqnarray}
-\bar D(Nx)\Delta^h &+& \big(\bar D(Nx)-\bar B(Nx)\big)\nabla_+^h\\
&\longrightarrow
&\frac{\sigma^2}{2}\frac{2(\alpha+1)(1-x)^2-(2\beta+1)x(2-x)}
{4(1-x)}\partI{}{x}
+\frac{\sigma^2}{2}\frac{x(2-x)}{4}\partII{}{x}\nonumber
\end{eqnarray}
in the limit $h=1/N\to0$.

Furthermore, the Racah polynomials, solutions to the finite
difference equation generated by (\ref{eq:DiscL}), converge to the
dual Jacobi polynomials up to some factor $\displaystyle
\frac{n!}{(\alpha+1)_n}$, since
\begin{eqnarray}
{_4F_3} \left(
\begin{matrix}
-n,  n+\alpha+\beta+1, -Nx, Nx-2N-1-\beta-\bar\beta \\
\alpha+1, -\beta-N, -N
\end{matrix}
\bigg\lvert 1 \right)\nonumber\\ \longrightarrow {_2F_1} \left(
\begin{matrix}
-n,  n+\alpha+\beta+1\\
\alpha+1
\end{matrix}
\bigg\lvert Z(x)\right)
\end{eqnarray}
as $N\to\infty$, with $Z(x)=x(2-x)$. We also have as $N\to\infty$,
\begin{equation}
\prod_{k=1}^{Nx}\frac{D(k)}{\bar D(k)} =
\prod_{k=1}^{Nx}\(1+\frac{\beta-\bar\beta}
{N(1-\frac{k}{N})+\bar\beta+1}\) \longrightarrow
\exp\(\int_0^x\frac{\beta-\bar\beta}{1-y}dy\)
=(1-x)^{\bar\beta-\beta}.
\end{equation}
Finally, from all the above limit relations, the Laplace transform
for the integral of the Racah process, given by (\ref{eq:RacahG})
with the appropriate assumptions, provides an extension with
underlying process $X_t$ on the lattice to the Laplace transform
(\ref{eq:dualJacobiG}) for the integral of the dual Jacobi process
in the particular case $\bar\alpha=\alpha$,
$\beta>\vert\bar\beta\vert$ and $\bar\beta\ge-1$.
\end{proof}


\section{Conclusion}
We have given a complete classification scheme for diffusion
processes for which Laplace transforms for integrals of stochastic
processes and transitional probability densities can be expressed as
integrals of hypergeometric functions against the spectral measure
for certain self-adjoint operators. The known models such as the
Ornstein-Uhlenbeck process, the CIR process and the geometric
Brownian motion fit into this classification scheme. We have also
presented extensions to these models in the quadratic
Ornstein-Uhlenbeck process and the Jacobi process. An extension of
the framework towards finite-state Markov processes related to
hypergeometric polynomials in the discrete series of the Askey
classification tree has been derived. Finally, we have explicitly
computed some limit relations between discrete and continuous
processes.

\bibliographystyle{plain}
\bibliography{stochint}

\end{document}